\documentclass[11pt,reqno,twoside]{amsart}

\usepackage{amsmath, amssymb, amscd, amsthm}
\usepackage{graphicx}

\theoremstyle{plain}
\newtheorem{thm}{Theorem}[section]

\newtheorem{prop}[thm]{Proposition}

\theoremstyle{remark}
\newtheorem{rem}[thm]{Remark}

\newcounter{sspar}[subsection]
\renewcommand\thesspar{(\thesubsection.\arabic{sspar})}
    {\par\ \newline
     \vskip-\baselineskip\vskip.1truecm
     \noindent\refstepcounter{sspar}
     \noindent\textbf{\thesspar} \ignorespaces}
    {\vskip-\baselineskip
    \ignorespaces}
    {\refstepcounter{sspar}
     \textup{\textbf{\thesspar}} \ignorespaces}
    {\vskip-\baselineskip
    \ignorespaces}
\newcommand{\R}{{\mathbb R}}

\newcommand{\C}{{\mathbb C}}
\newcommand{\Z}{{\mathbb Z}}
\newcommand{\al}{\alpha}

\newcommand{\de}{\delta}
\newcommand{\De}{\Delta}
\newcommand{\eps}{\varepsilon}
\newcommand{\si}{\sigma}

\newcommand{\te}{\theta}
\newcommand{\Te}{\Theta}
\newcommand{\la}{\lambda}
\newcommand{\La}{\Lambda}

\newcommand{\Ups}{\Upsilon}

\newcommand{\Om}{\Omega}
\newcommand{\lt}{\ell^2}
\newcommand{\inprod}[2]{\langle #1,#2 \rangle}
\newcommand{\tensor}{\otimes}
\newcommand{\dirint}{\sideset{}{^\oplus}\int\limits_{0}^{-\frac{\pi}{2\ln
q}}}

\newcommand{\ph}[5]{\,_{#1}\varphi_{#2}\!\left( \genfrac{.}{.}{0pt}{}{#3}{#4}
\,;#5 \right)}
\newcommand{\p}{\tilde{p}}
\newcommand{\Pt}{\widetilde{P}}
\newcommand{\hf}{\frac{1}{2}}
\newcommand{\su}{\mathfrak{su}}
\newcommand{\s}{\tilde{s}}
\newcommand{\w}{\overline{w}}

\newcommand{\M}{\hat{M}}
\DeclareMathOperator{\sgn}{\mathrm{sgn}}
\numberwithin{equation}{section}

\begin{document}

\title{Bilinear summation formulas from quantum algebra representations}
\author{Wolter Groenevelt}
\date{January 28, 2002}
\address{Technische Universiteit Delft, ITS-TWA \\
Postbus 5031, 2600 GA Delft, The Netherlands}
\email{W.G.M.Groenevelt@its.tudelft.nl}

\begin{abstract}
The tensor product of a positive and a negative discrete series representation of the quantum algebra $U_q\big(\su(1,1)\big)$ decomposes as a direct integral over the principal unitary series representations. Discrete terms can appear, and these terms are a finite number of discrete series representations, or one complementary series representation. From the interpretation as overlap coefficients of little $q$-Jacobi functions and Al-Salam and Chihara polynomials in base $q$ and base $q^{-1}$, two closely related bilinear summation formulas for the Al-Salam and Chihara polynomials are derived. The formulas involve Askey-Wilson polynomials, continuous dual $q$-Hahn polynomials and little $q$-Jacobi functions. The realization of the discrete series as $q$-difference operators on the spaces of holomorphic and anti-holomorphic functions, leads to a bilinear generating function for a certain type of $_2\varphi_1$-series, which can be considered as a special case of the dual transmutation kernel for little $q$-Jacobi functions.
\end{abstract}
\maketitle

\section{Introduction}
The Askey-scheme of hypergeometric and basic hypergeometric orthogonal polynomials, see \cite{KS}, consists of polynomial systems which can be defined in terms of hypergeometric or basic hypergeometric functions. On top of the Askey-scheme is a four parameter family of orthogonal polynomials, introduced by Askey and Wilson in \cite{AW}, called the Askey-Wilson polynomials. The other families can be derived by limit transitions from the Askey-Wilson polynomials. The Al-Salam and Chihara polynomials are Askey-Wilson polynomials with two parameters equal to zero. In this paper we use representation theory of the quantized universal enveloping algebra $U_q\big(\su(1,1)\big)$ to derive two bilinear summation formulas for the Al-Salam and Chihara polynomials in base $q$ and $q^{-1}$, involving Askey-Wilson polynomials, continuous dual $q$-Hahn polynomials and little $q$-Jacobi functions, and we discuss some consequences.

Most of the hypergeometric polynomials in the Askey-scheme are related to the representation theory of Lie groups and Lie algebras, e.g. the polynomials appear as matrix coefficients of irreducible unitary representations, see e.g. \cite{Koo} and \cite{VK} and references therein. In the same way most of the basic hypergeometric polynomials are related to representation theory of quantum groups and quantum algebras. Using representation theory one can obtain several identities, such as generating functions or convolution identities, for special functions of (basic) hypergeometric type. One method to find such identities uses representations of a Lie algebra, and is due to Granovskii and Zhedanov \cite{GZ}. The idea is to consider (generalized) eigenvectors of a certain element of the Lie algebra, which acts as a recurrence operator in an irreducible representation. From the tensor product decomposition, one can find identities for the eigenvectors. Especially one finds identities for the special functions that appear as overlap coefficients. This idea is used by Koelink and Van der Jeugt \cite{KJ}, where tensor products of positive discrete series representation of the Lie algebra $\su(1,1)$ are used to obtain convolution identities for orthogonal polynomials. In \cite{KJ} the method is also applied to tensor products of positive discrete series representations of the quantum algebra $U_q\big(\su(1,1)\big)$. In that case the operators considered, are related to so-called twisted primitive elements in order to control the action in the tensor product representations.

In \cite{GK} the method is used on the tensor product of a positive and a negative discrete series representation of the Lie algebra $\su(1,1)$. In that case non-polynomial hypergeometric functions are needed, namely the Meixner function as defined in \cite{MR}, see also \cite{Koe}. These Meixner functions can be considered as non-polynomial extensions of the Meixner polynomials. The goal of this paper is to find $q$-analogues of the results of \cite{GK} using representation theory of $U_q\big(\su(1,1)\big)$. The method we use is different from the method applied in \cite{GK}.

The irreducible unitary representations of $U_q\big(\su(1,1)\big)$ are the discrete series representations, which act on $\ell^2(\Z_{\geq 0})$, and the principal unitary series, the complementary series and the strange series representations, which all act on $\ell^2(\Z)$. We consider the tensor product of a positive and a negative discrete series representation. This tensor product decomposes as a direct integral over the unitary representations, see Kalnins and Miller \cite{KM}. A finite number of discrete terms can appear in the decomposition, and these discrete terms are discrete series or at most one complementary series representation. The strange series do not appear in the decomposition. As overlap coefficients related to the positive and negative discrete series representations, we find Al-Salam and Chihara polynomials in base $q^2$, respectively in base $q^{-2}$. For the principal unitary series, we find little $q$-Jacobi functions, see \cite{Ka}, \cite{KSt1}. So in this context, the little $q$-Jacobi functions can be considered as non-polynomial extensions of the Al-Salam and Chihara polynomials, see also \cite{KSt2}.

In section \ref{sec2} we give the definition of the Askey-Wilson polynomials, their orthogonality relations and we give an (apparently new) generating function for these polynomials. This generating function plays a key role in this paper. Then we determine the exact decomposition of the tensor product of a positive and a negative discrete series representation of $U_q\big(\su(1,1)\big)$ by considering the action of the Casimir element. We find continuous dual $q$-Hahn polynomials as Clebsch-Gordan coefficients. In section \ref{sec3} we determine (generalized) eigenvectors of a certain element in the various irreducible representations and in the tensor product representation. We determine the Clebsch-Gordan coefficients for the bases of eigenvectors, which turn out to be Askey-Wilson polynomials. As a result, see theorems \ref{thm1} and \ref{thm2}, we obtain two summation formulas involving Al-Salam and Chihara polynomials, continuous dual $q$-Hahn polynomials, Askey-Wilson polynomials and little $q$-Jacobi functions. In section \ref{sec4} we realize the generators of $U_q\big(\su(1,1)\big)$ in the positive and negative discrete series representations as $q$-difference operators on the space of holomorphic, respectively anti-holomorphic functions. The eigenvectors for the discrete series now become known generating functions for Al-Salam and Chihara polynomials. From these realization of the eigenvectors, we derive a bilinear generating function for a certain type of $_2\varphi_1$-series, see theorem \ref{thm4.4} This gives a quantum group theoretical proof of a special case of the dual transmutation kernel for the little $q$-Jacobi functions, which has recently been found by Koelink and Rosengren \cite{KR}.\\

\emph{Notations.}
Throughout this paper we assume $0<q<1$. We use the notation for basic hypergeometric series and $q$-shifted factorials as in the book of Gasper and Rahman \cite{GR}, i.e.
\[
\begin{split}
\ph{r}{s}{a_1,a_2,\ldots,a_r}{b_1,\ldots,b_s}{q,z} &= \sum_{k=0}^\infty
\frac{ (a_1,\ldots,a_r;q)_k }{ (q,b_1,\ldots, b_s;q)_k} \left( (-1)^k
q^{\hf k(k-1)} \right)^{1+s-r}z^k, \\ (a_1,\ldots,a_r;q)_k &= (a_1;q)_k
 \ldots (a_r;q)_k, \quad (a;q)_k = \prod_{j=0}^{k-1} (1-aq^j).
\end{split}
\]
The $_r\varphi_s$-series converges absolutely for all $z$ if $r\leq s$, for $|z|<1$ if $r=s+1$ and diverges for $r>s+1$. The $_{s+1}\varphi_s$-series has a unique analytic continuation to $\C\setminus[1,\infty)$, see \cite[\S4.5]{GR}. We often use the analytic continuation implicitly.

A basic hypergeometric series is called very-well-poised if $r=s+1$ and
$a_1q=a_4b_3 = a_5b_4=\ldots =a_{s+1}b_s$, $a_2=q\sqrt{a_1}$ and
$a_3=-q\sqrt{a_1}$. We use the notation for a very-well-poised basic hypergeometric series
as in \cite[\S2.1]{GR}
\[
\begin{split}
_{s+1}W_s(a_1;a_4,\ldots,a_{s+1};q,z) &= \ph{s+1}{s}{a_1, q\sqrt{a_1}, -q\sqrt{a_1},
a_4, \ldots, a_{s+2} }{ \sqrt{a_1}, -\sqrt{a_1}, a_1q/a_4, \ldots, a_1q/a_{s+1}}
{q,z} \\&= \sum_{k=0}^\infty \frac{ 1-a_1q^{2k} }{1-a_1} \frac{ (a_1,a_4, \ldots, a_{s+1};q)_k \,z^k}
{ (q, a_1q/a_4, \ldots, a_1q/a_{s+1} ;q)_k } .
\end{split}
\]

If $dm$ is a positive measure on $\R$, we denote by $dm^\hf$ the positive measure with the property that $dm$ is the product measure of $dm^\hf$ with itself, restricted to the diagonal.

\emph{Acknowledment.} I thank Erik Koelink for comments on previous versions.

\section{Decomposition of tensor product representations} \label{sec2}
In this section we consider the quantized universal enveloping
algebra $U_q\big(\su(1,1)\big)$ and the irreducible
representations. We decompose the tensor product of a positive and a negative discrete series representation into a direct integral of principal unitary series. Under certain conditions discrete terms appear in the decomposition. These discrete terms are a finite number of negative discrete series representations, or one complementary series representation. Also we find continuous dual $q$-Hahn polynomials as Clebsch-Gordan coefficients, cf. Kalnins and Miller \cite{KM}.

\subsection{The quantized universal enveloping
algebra $\boldsymbol{U_q\big(\su(1,1)\big)}$}
$U_q\big(\su(1,1)\big)$ is the unital, associative, complex algebra generated by $A$, $B$, $C$ and $D$, subject to the relations
\[
AD = 1 = DA, \quad AB = qBA, \quad AC= q^{-1}CA, \quad BC-CB =
\frac{A^2-D^2}{q-q^{-1}}.
\]
The Casimir element
\begin{equation} \label{Casimir}
\Om = \frac{q^{-1} A^2 +qD^2-2}{(q^{-1}-q)^2} + BC =
\frac{q^{-1}D^2+qA^2-2}{(q^{-1}-q)^2} +CB
\end{equation}
is a central element of $U_q\big(\mathfrak{su}(1,1)\big)$. The algebra
$U_q\big(\mathfrak{su}(1,1)\big)$ is a Hopf $*$-algebra with
comultiplication $\De$ given by
\begin{equation} \label{comult}
\De(A) = A \tensor A, \quad \De(B) = A \tensor B + B \tensor D,
\quad \De(C) = A \tensor C + C \tensor D, \quad \De(D) = D \tensor
D.
\end{equation}
The $*$-structure is defined by
\[
A^*=A, \quad B^*=-C, \quad C^* = -B, \quad D^* = D.
\]
The irreducible unitary representations have been determined by Burban
and Klimyk \cite{BK}. There are five classes of irreducible unitary
representations of $U_q\big(\mathfrak{su}(1,1)\big)$:\\

\emph{Positive discrete series.} The positive discrete series $\pi^+_k$
are labeled by $k>0$. The representation space is $\lt(\Z_{\geq 0})$
with orthonormal basis $\{e_n\}_{n \in \Z_{\geq 0} }$. The action is
given by
\begin{equation} \label{pos}
\begin{split}
\pi^+_k(A) e_n &= q^{k+n} e_n, \quad \pi^+_k(D) e_n =
q^{-(k+n)}e_n,
 \\ (q^{-1}-q)\pi^+_k(B) e_n &= q^{-\hf-k-n}
\sqrt{(1-q^{2n+2})(1-q^{4k+2n})}\, e_{n+1},  \\
(q^{-1}-q)\pi^+_k(C) e_n &= -q^{\hf-k-n}
\sqrt{(1-q^{2n})(1-q^{4k+2n-2})}\, e_{n-1},
\\ (q^{-1}-q)^2 \pi^+_k(\Om) e_n &= (q^{2k-1}+q^{-(2k+1)}-2)e_n.
\end{split}
\end{equation}

\emph{Negative discrete series.} The negative discrete series
$\pi^-_k$ are labeled by $k>0$. The representation space is
$\lt(\Z_{\geq 0})$ with orthonormal basis $\{e_n\}_{n \in \Z_{\geq
0}}$. The action is given by
\begin{equation} \label{neg}
\begin{split}
\pi^-_k(A) e_n &= q^{-(k+n)} e_n, \quad \pi^-_k(D) e_n =
q^{k+n}e_n,
 \\  (q^{-1}-q)\pi^-_k(B) e_n &= -q^{\hf-k-n}
\sqrt{(1-q^{2n})(1-q^{4k+2n-2})}\, e_{n-1},
\\(q^{-1}-q)\pi^-_k(C) e_n &= q^{-\hf-k-n}
\sqrt{(1-q^{2n+2})(1-q^{4k+2n})}\, e_{n+1}, \\ (q^{-1}-q)^2
\pi^-_k(\Om) e_n &= (q^{2k-1}+q^{-(2k+1)}-2)e_n.
\end{split}
\end{equation}

\emph{Principal unitary series.} The principal unitary series
representations $\pi^P_{\rho,\eps}$ are labeled by $0\leq\rho \leq
-\frac{\pi}{2 \ln q}$ and $\eps \in [0,1)$, where $(\rho,\eps)
\neq (0,\hf)$. The representation space is $\lt(\Z)$ with
orthonormal basis $\{e_n\}_{n \in \Z}$. The action is given by
\begin{equation} \label{prin}
\begin{split}
\pi^P_{\rho,\eps}(A) e_n &= q^{n+\eps} e_n, \quad
\pi^P_{\rho,\eps} (D) e_n = q^{-(n+\eps)}e_n,
\\  (q^{-1}-q)\pi^P_{\rho,\eps}(B) e_n
&= q^{-\hf-n-\eps} \sqrt{(1-q^{2n+2\eps+2i\rho+1})
(1-q^{2n+2\eps-2i\rho+1})}\, e_{n+1},
\\(q^{-1}-q)\pi^P_{\rho,\eps}(C) e_n &= -q^{\hf-n-\eps}
\sqrt{(1-q^{2n+2\eps+2i\rho-1}) (1-q^{2n+2\eps-2i\rho-1})}\, e_{n-1},\\
(q^{-1}-q)^2 \pi^P_{\rho,\eps} (\Om) e_n &=
(q^{2i\rho}+q^{-2i\rho} -2)e_n.
\end{split}
\end{equation}
For $(\rho,\eps)= (0,\hf)$ the representation $\pi^P_{0,\hf}$ splits into a direct
sum of a positive and a negative discrete series representation: $\pi^P_{0,\hf}= \pi^+_\hf
\oplus \pi^-_\hf$. The representation space splits into two invariant subspaces: $\{e_n \, |\,
n<0 \} \oplus \{ e_n\, | \, n \geq 0 \}$.

\emph{Complementary series.} The complementary series
representations $\pi^C_{\la,\eps}$ are labeled by $\la$ and
$\eps$, where $\eps \in [0,\hf)$ and $\la \in (-\hf,-\eps)$, or
$\eps \in (\hf,1)$ and $\la \in (-\hf, \eps-1)$. The
representation space is $\lt(\Z)$ with orthonormal basis
$\{e_n\}_{n \in \Z}$. The action is given by
\begin{equation} \label{comp}
\begin{split}
\pi^C_{\la,\eps}(A) e_n &= q^{n+\eps} e_n, \quad \pi^C_{\la,\eps}
(D) e_n = q^{-(n+\eps)}e_n,
\\  (q^{-1}-q)\pi^C_{\la,\eps}(B) e_n
&= q^{-\hf-\eps-n}\sqrt{(1-q^{2n+2\eps+2\la+2})
(1-q^{2n+2\eps-2\la})}\, e_{n+1},
\\(q^{-1}-q)\pi^C_{\la,\eps}(C) e_n &= -q^{\hf-\eps-n} \sqrt{
(1-q^{2n+2\eps+2\la}) (1-q^{2n+2\eps-2\la-2})}\, e_{n-1},\\
(q^{-1}-q)^2 \pi^C_{\la,\eps} (\Om) e_n &=
(q^{2\la+1}+q^{-(2\la+1)}-2)e_n.
\end{split}
\end{equation}

The fifth class consists of the strange series representations. The
strange series representations do not appear in the decomposition of
the tensor product of a positive and a negative series representation,
therefore we do not need them in this paper.

Note that the operators are unbounded, with common domain the set of finite
linear combinations of the basisvectors. The operators in \eqref{pos}-\eqref{comp} define $*$-representations in the sense of Schm\"udgen \cite[Ch.8]{Sch}.

\subsection{Tensor product of positive and negative discrete
series representations} The decomposition of the tensor product of a
positive and a negative discrete series representation, has been
determined by Kalnins and Miller in \cite{KM}. They find continuous
dual $q$-Hahn polynomials as Clebsch-Gordan coefficients. The
continuous dual $q$-Hahn polynomials are a subclass of the Askey-Wilson
polynomials $p_n$, defined by (see Askey and Wilson \cite{AW}, Koekoek
and Swarttouw \cite{KS})
\begin{equation} \label{AW pol}
p_n(\cos \te;a,b,c,d|q)= a^{-n} (ab,ac,ad;q)_n \ph{4}{3}{q^{-n},
abcdq^{n-1}, ae^{i\te}, ae^{-i\te} } {ab, ac, ad} {q, q} .
\end{equation}
By Sear's $_4\varphi_3$ transformation formula \cite[eq.(III.16)]{GR}
the polynomials $p_n$ are symmetric in the parameters $a$, $b$, $c$ and
$d$. Let $a,b,c,d$ be real, or appearing in complex conjungate pairs,
and let the pairwise products of $a,b,c,d$ be smaller than $1$, then
the Askey-Wilson polynomials are orthogonal with respect to a positive
measure supported on a subset of $\R$. The orthonormal Askey-Wilson
polynomials $\p_n$ are defined by
\begin{equation} \label{on AW}
\p_n(y;a,b,c,d|q) = \sqrt{
\frac{ (abcd;q)_{2n} } { (q,ab,ac,ad,bc,bd,cd,abcdq^{n-1};q)_n }}
p_n(y;a,b,c,d|q).
\end{equation}
They are orthonormal with respect to the measure $dm(\cdot;a,b,c,d|q)$
given by
\begin{gather}
\int_\R f(y) dm(y;a,b,c,c|q) =\int_{0}^\pi f(\cos\te) w(\cos\te) d\te +
\sum_k f(x_k) w_k, \label{measure}
\\ w(\cos\te)=w(\cos\te;a,b,c,d|q) = \frac{1}{2\pi} \frac{
(q,ab,ac,ad,bc,bd,cd;q)_\infty} { (abcd;q)_\infty } \left| \frac{
(e^{2i\te};q)_\infty } { (ae^{i\te}, be^{i\te}, ce^{i\te},
de^{i\te};q)_\infty } \right|^2, \nonumber
\end{gather}
where $x_k=\mu(eq^k)$ for $e$ any of the parameters $a$, $b$, $c$, $d$.
Here and elsewhere $\mu(y) = \hf(y+y^{-1})$.
The sum is over $k \in \Z_{\geq 0}$ such that $|eq^k|>1$. If we
assume $e=a$, we have
\[
w_k = w_k(a;b,c,d|q) =
\frac{1-a^2q^{2k}}{1-a^2} \frac{ (a^{-2}, bc, bd, cd
;q)_\infty }{(b/a, c/a, d/a, abcd;q)_\infty} \frac{(a^2,ab,ac,ad;q)_k}
{(q, aq/b, aq/c, aq/d ; q)_k} \left(\frac{1}{abcd}\right)^k.
\]
For future references we give an (apparently new) generating function for the Askey-Wilson
polynomials.
\begin{thm} \label{thm2.1}
For $|t|<1$ the Askey-Wilson polynomials satisfy the following
generating function
\begin{equation} \label{gen AW}
\begin{split}
\sum_{n=0}^\infty &\frac{ (abcd;q)_{2n} p_n(\cos\te;a,b,c,d|q)}
{(q,ab,ac,bc,abcdq^{n-1};q)_n } \frac{ (r/t, abc/r;q)_n } {
(abcdt/r, rd;q)_n} t^n =\\&  \frac{ (abcd, dt, abcte^{i\te}/r, re^{i\te} ;q)_\infty } { (abcdt/r, dr, abce^{i\te}, te^{i\te};q)_\infty }\
_8W_7(abce^{i\te}/q; ae^{i\te}, be^{i\te}, ce^{i\te}, r/t, abc/r ;q, te^{-i\te}).
\end{split}
\end{equation}
\end{thm}
\begin{proof}
We start with the sum $S$ on the left hand side of \eqref{gen AW}. Using the asymptotic behaviour of the Askey-Wilson polynomials, see \cite[eq.(7.5.13)]{GR}, we find that $S$ converges absolutely for $|t|<1$. We use \eqref{AW pol} and Watson's transformation \cite[eq.(III.19)]{GR}, to write the Askey-Wilson polynomial as a multiple of a very-well-poised $_8\varphi_7$-series;
\begin{multline*}
p_n(\cos\te;a,b,c,d|q) =\\ \frac{(ab, ac, bc, de^{-i\te};q)_n } { (abce^{i\te};q)_n } e^{in\te}\,
_8W_7(abce^{i\te}/q; ae^{i\te}, be^{i\te}, ce^{i\te}, abcdq^{n-1}, q^{-n};q, qe^{-i\te}/d).
\end{multline*}
Next we write out the $_8\varphi_7$-series as a sum, so $S$ becomes a double sum, which is absolutely convergent. We interchange summations
\[
\sum_{n=0}^\infty \sum_{l=0}^n = \sum_{l=0}^\infty \sum_{p=0}^\infty \ , \qquad p=n-l,
\]
and we use
\begin{align*}
\frac{ (q^{-l-p};q)_l }{ (q;q)_{l+p} } &= (-1)^l \frac{ q^{-\hf l(l+1)-lp} }{ (q;q)_p }, \\
\frac{ (\al;q)_{2l+2p}\, (\al q^{l+p-1};q)_l } {(\al q^{l+p-1};q)_{l+p} } &= \frac{ 1-\al q^{2l+2p-1} }
{ 1- \al q^{2l-1} } (\al;q)_{2l} \, (\al q^{2l-1};q)_p, \\
(q^{1-l-p}/\al ;q)_l &= (-1)^l q^{-\hf l(l-1)-lp} \al^{-l} \frac{ (\al q^l;q)_p \, (\al;q)_l }
{ (\al;q)_p }, \\
(\al q^{p} ;q)_l &= \frac{ (\al q^l;q)_p\, (\al;q)_l }{ (\al;q)_p },
\end{align*}
then after some cancellations we have
\[
S = \sum_{l=0}^\infty  \frac{ 1-abcq^{2l-1}e^{i\te} } { 1-abce^{i\te}/q} \frac{ (abcd;q)_{2l}\,
(ae^{i\te}, be^{i\te}, ce^{i\te}, abce^{i\te}, abc/r, r/t ;q)_l } { (q, ab, ac, bc, abce^{i\te},
abcq^le^{i\te}, abcdt/r ;q)_l } e^{-il\te} t^l S_l,
\]
where $S_l$ is the sum over $p$. We write $S_l$ as a very-well-poised $_6\varphi_5$-series, which is
summable by Jackson's summation formula \cite[eq.(II.20)]{GR};
\[
S_l ={} _6W_5(abcdq^{2l-1}; rq^l/t, abcq^l/r, de^{-i\te}; q, e^{i\te}t) = \frac{ (abcdq^{2l}, dt, abctq^l e^{i\te}/r,
 rq^l e^{i\te} ;q)_\infty } { (abcdtq^l/r, rdq^l, abcq^{2l}e^{i\te}, te^{i\te};q)_\infty}.
\]
Now $S$ reduces to a single sum, which turns out to be a multiple of a very-well-poised $_8\varphi_7$-series;
\[
S= \frac{ (abcd, dt, abcte^{i\te}/r, re^{i\te} ;q)_\infty } { (abcdt/r, dr, abce^{i\te}, te^{i\te};q)_\infty }{}
_8W_7(abce^{i\te}/q; ae^{i\te}, be^{i\te}, ce^{i\te}, r/t, abc/r ;q, te^{-i\te}).
\]
This is the desired result.
\end{proof}
\begin{rem} \label{rem2.2}
The $_8\varphi_7$-series on the right hand side of \eqref{gen AW} can be written as the sum of two balanced $_4\varphi_3$-series by \cite[eq.(III.36)]{GR}. For $t=q/\al$ and $r=q^{1-m}/\al$ this reduces to one balanced $_4\varphi_3$-series, and this gives the well known connection formula, see \cite[\S6]{AW}, \cite[\S7.6]{GR},
\[
\begin{split}
&p_m(\cos\te;a,b,c,\al|q) = \sum_{n=0}^m c_{n,m} p_n(\cos\te;a,b,c,d|q), \\&
c_{n,m} = \frac{ (q^{-m},abc\al q^{m-1};q)_n} { (q, abcdq^{n-1};q)_n } \frac{ (abq^n, acq^n, bcq^n, \al/d;q)_{m-n} } { (abcdq^{2n};q)_{n-m} } (-1)^n d^{m-n} q^{mn-\hf n(n-1)}.
\end{split}
\]
\end{rem}

The continuous dual $q$-Hahn polynomials $P_n$ are obtained from the
Askey-Wilson polynomials by taking $d=0$;
\begin{equation} \label{CDH}
 P_n(\cos\te)=P_n(\cos\te;a,b,c|q) =
p_n(\cos\te;a,b,c,0|q)= a^{-n} (ab,ac;q)_n \ph{3}{2}{q^{-n}, ae^{i\te},
ae^{-i\te}} {ab,ac} {q,q}.
\end{equation}
The orthonormal continuous dual $q$-Hahn polynomials $\Pt_n$ are
orthonormal with respect to the measure
$dm(\cdot;a,b,c|q)=dm(\cdot;a,b,c,0|q)$ and they satisfy the recurrence
relation
\begin{equation} \label{rec cdH}
2y\Pt_n(y) = a_n \Pt_{n+1}(y) + b_n \Pt_n(y) + a_{n-1} \Pt_{n-1}(y),
\end{equation}
where
\begin{align*}
a_n & = \sqrt{ (1-q^{n+1})(1-abq^{n})(1-acq^n)(1-bcq^n) }, \\ b_n &=
a+a^{-1}-a^{-1}(1-abq^n)(1-acq^n) - a(1-q^n)(1-bcq^{n-1}).
\end{align*}

In \cite{KM} the decomposition of the tensor product $\pi^+_{k_1}
\tensor \pi^-_{k_2}$ is found by considering the action of the
Casimir element $\Om$ in this tensor product representation. Since
our main focus is on special functions, we need to know the
Clebsch-Gordan decomposition and the matrix elements of the
intertwiner exactly. Therefore we repeat the proof given in
\cite{KM} in somewhat more detail.

From \eqref{Casimir} and \eqref{comult} we find
\[
\begin{split}
\De(\Om) = &\frac{1}{(q^{-1}-q)^{2}} \big[ q^{-1}( A^2 \tensor A^2 )+ q
(D^2 \tensor D^2) -2 (1\tensor 1) \big] \\&+ A^2 \tensor BC +AC\tensor
BD + BA \tensor DC + BC \tensor D^2.
\end{split}
\]
We define elements $f_n^p \in \lt(\Z_{\geq 0}) \tensor
\lt(\Z_{\geq 0})$ by
\[
f_n^p =
\begin{cases} e_n \tensor e_{n-p}, & p \leq 0,\\ e_{n+p} \tensor
e_n, & p \geq 0,
\end{cases}
\]
and we define the space $H_p$ by
\[
H_p = \overline{ \C\{f_n^p|n \in \Z_{\geq 0}\}} \cong \lt(\Z_{\geq
0}).
\]
For fixed $p$ we let $\De(\Om)$ act on finite linear combinations
of elements $f_n^p$. We see by a straightforward computation that
$(q^{-1}-q)^2 \pi^+_{k_1} \tensor \pi^-_{k_2}\big( \De(\Om)
\big)+2$ can be identified with the three term recurrence relation
for the continuous dual $q$-Hahn polynomials \eqref{rec cdH} in
base $q^2$ with parameters
\begin{equation} \label{abc}
a  =
\begin{cases} q^{2k_2-2k_1-2p+1},&  p \leq 0,\\
q^{2k_1-2k_2+2p+1},& p \geq 0,
\end{cases}
\quad b  = q^{2k_1+2k_2-1}, \quad c =
\begin{cases} q^{2k_1-2k_2+1}, & p \leq 0, \\
q^{2k_2-2k_1+1}, & p \geq 0.
\end{cases}
\end{equation}
\begin{prop} \label{prop2.1}
The operator $\La_p$ defined by
\[
\begin{split}
\La_p : H_p &\rightarrow L^2\big(\R,dm(\cdot;a,b,c|q^2)\big) \\
f_n^p &\mapsto \Pt_n(\cdot;a,b,c|q^2)
\end{split}
\]
is unitary and intertwines $\pi^+_{k_1} \tensor \pi^-_{k_2}\big(
\De(\Om) \big)$ acting on $H_p$ with $(q^{-1}-q)^{-2}M_{2x-2}$
acting on \\ $L^2\big(\R,dm(x;a,b,c|q^2)\big)$ with $a$, $b$, $c$ defined in
\eqref{abc}.
\end{prop}
Here and elsewhere $M$ denotes the multiplication operator, i.e.
$M_f g(x) = f(x)g(x)$.
\begin{proof}
$\pi^+_{k_1} \tensor \pi^-_{k_2}\big( \De(\Om) \big)$ resticted to
$H_p$ extends to a bounded self-adjoint Jacobi operator on $H_p$,
see Akhiezer \cite{Akh}. The intertwining now follows from
\eqref{rec cdH}. Since $\La_p$ maps an orthonormal basis onto
another, $\La_p$ is unitary.
\end{proof}

We define a map $\vartheta: U_q\big(\su(1,1)\big) \rightarrow U_q\big(\su(1,1)\big)$ by
\[
\vartheta(A)=D, \qquad \vartheta(B)=C, \qquad \vartheta(C)=B, \qquad \vartheta(D)=A.
\]
Then $\vartheta$ is an algebra homomorphism and from \eqref{pos} and \eqref{neg} we find $\pi^\pm_k\big(\vartheta(X)\big) = \pi^\mp_k(X)$ for $X \in U_q\big(\su(1,1)\big)$. From \eqref{comult} follows that $\vartheta$ is an anti-coalgebra homomorphism, i.e. $\De \circ \vartheta = \si \circ \De$, where $\si: U_q\big(\su(1,1)\big) \rightarrow U_q\big(\su(1,1)\big)$ denotes the flip automorphism, $\si(u\tensor v) = v \tensor u$. Since $U_q\big(\su(1,1)\big)$ is a quasitriangular Hopf-algebra we have $\pi^+_{k_1} \tensor \pi^-_{k_2} \cong \pi^-_{k_2} \tensor \pi^+_{k_1}$, where the intertwiner is induced by the universal $R$-matrix, see \cite[\S4.2, \S6.4]{CP}. So $(\pi^+_{k_1} \tensor \pi^-_{k_2})\circ \De \circ \vartheta$ is equivalent to the standard tensor product representation $(\pi^+_{k_2} \tensor \pi^-_{k_1}) \circ \De$. This shows that the case $k_2 \geq k_1$ gives results similar to the case $k_1 \geq k_2$. So from here on we assume $k_2 \geq k_1$.

From proposition \ref{prop2.1} follows that the spectrum of $\pi^+_{k_1} \tensor \pi^-_{k_2}\big( \De(\Om)\big)$ can be read off from the support of the orthogonality measure $dm(\cdot;a,b,c|q^2)$, with $a$, $b$, $c$ as in \eqref{abc}. The measure always has an absolutely continuous part, and possibly a finite set of discrete mass points when one of the parameters is greater than one. We distinguish three different cases.
\begin{enumerate}
\item If $k_1-k_2 \geq -\frac{1}{2}$ and $k_1+k_2 \geq \frac{1}{2}$ the
measure $dm$ in proposition \ref{prop2.1} is absolutely continuous
and has support $[-1,1]$ for all $p \in \Z$. For this part we
recognize the action of $\Om$ in the principal unitary series from
\eqref{prin}, using $e^{i\te}=q^{2i\rho}$. From the action of $A$
in the tensor product representation, $\pi^+_{k_1} \tensor \pi^-_{k_2}\big( \De(A)\big) = \pi^+_{k_1}(A) \tensor \pi^-_{k_2}(A)$, we find $\eps=k_1-k_2+L$,
where $L$ is the unique non-negative integer such that $\eps \in
[0,1)$.

\item If $k_1+k_2<\hf$ the measure $dm$ in proposition \ref{prop2.1} has one discrete mass point
outside $[-1,1]$ for all $p \in \Z$. In this case
$(q^{2k_1+2k_2-1}+q^{1-2k_1-2k_2}-2)/(q^{-1}-q)^2$ is also an
eigenvalue of $\pi^+_{k_1} \tensor \pi^-_{k_2}\big( \De(\Om) \big)$. We
recognize the action of $\Om$ in the complementary series
representation from \eqref{comp} with $\la = -k_1-k_2$. Again from the action of $A$ in the tensor product representation we
find $\eps=k_1-k_2+L$.

\item If $k_1-k_2 < - \hf$, then the support of $dm$ in proposition \eqref{prop2.1} contains finitely
many points outside $[-1,1]$. The number of discrete points depends on $p$. These discrete mass points correspond to
eigenvalues of $\pi^+_{k_1} \tensor \pi^-_{k_2}\big( \De(\Om) \big)$ of
the form $(q^{2k_1-2k_2+1+2j)} + q^{-2k_1+2k_2-1-2j)}+2)/(q^{-1}-q)^2$.
Here $j=0,\ldots, K$ for $p \leq 0$, and $K$ is the largest integer
such that $k_1-k_2+\hf+K<0$. For $0 \leq p \leq K$, we have
$j=p,\ldots,K$ and for $p >K$ there are no discrete mass points. Here we
recognize the action of $\Om$ in a discrete series representation from \eqref{pos} and \eqref{neg}. From
the action of $A$ in the tensor product representation we find that this is a negative discrete series
representation with label $k_2-k_1-j$.
\end{enumerate}
We have the following decomposition.
\begin{thm} \label{decomp}
For $k_1 \leq k_2$ the decomposition of the tensor product of positive
and negative discrete series representations of $U_q\big(\su(1,1)\big)$
is
\begin{align*}
\pi_{k_1}^+ \tensor\pi_{k_2}^- &\cong \dirint \pi^P_{\rho,\eps} d \rho,
& k_1-k_2 \geq -\frac{1}{2}, k_1+k_2 \geq \frac{1}{2}, \\ \pi_{k_1}^+
\tensor\pi_{k_2}^- &\cong \dirint \pi_{\rho,\eps}^P d \rho \oplus
\pi^C_{\la, \eps}, & k_1+k_2<\frac{1}{2}, \\ \pi_{k_1}^+
\tensor\pi_{k_2}^- &\cong \dirint \pi^P_{\rho,\eps} d \rho \oplus
\bigoplus_{\substack{j \in \Z_{\geq 0}\\ k_2-k_1-\frac{1}{2}-j>0}}
\pi_{k_2-k_1-j}^-, & k_1-k_2 < -\frac{1}{2},
\end{align*}
where $\eps = k_1-k_2+L$, $L$ is the unique integer such that $\eps \in
[0,1)$ and $\la = -k_1-k_2$. Further, under the identification above,
for $y=\hf(q^{2i\rho}+q^{-2i\rho})$,
\begin{equation} \label{CG decomp}
e_{n_1} \tensor e_{n_2} =
\begin{cases}
\displaystyle (-1)^{n_1-n_2}\int_\R \Pt_{n_1}(y;a,b,c|q^2)
e_{n_1-n_2-L}dm^\hf(y;a,b,c|q^2)
 , & n_1 \leq n_2, \\ \displaystyle
\int_\R \Pt_{n_2}(y;a,b,c|q^2) e_{n_1-n_2-L}dm^\hf(y;a,b,c|q^2)
 , & n_1 \geq n_2,
\end{cases}
\end{equation}
where $\Pt_n$ is an orthonormal continuous dual $q$-Hahn
polynomial with parameters $a,b,c$ given by
\begin{equation} \label{param}
a  =
\begin{cases} q^{2k_2-2k_1+2n_2-2n_1+1},&  n_1 \leq n_2,\\
q^{2k_1-2k_2+2n_1-2n_2+1},& n_1 \geq n_2,
\end{cases}
\quad b  = q^{2k_1+2k_2-1}, \quad c =
\begin{cases} q^{2k_1-2k_2+1}, & n_1 \leq n_2, \\
q^{2k_2-2k_1+1}, & n_1 \geq n_2.
\end{cases}
\end{equation}
\end{thm}
We can also give the inverse of \eqref{CG decomp} explicitly, e.g.
if no discrete terms occur in the decomposition of the tensor
product, the decomposition of
\[
f \tensor e_{p-L} = \int_{-1}^1 f(y) e_{p-L} dy\ \in\
L^2(-1,1)\tensor \lt(\Z) \cong \dirint \lt(\Z) d\rho,
\]
is given by
\begin{equation} \label{inv CG}
f \tensor e_{p-L} =
\begin{cases}
\displaystyle (-1)^p \sum_{n=0}^\infty \left[\int_{-1}^1
\Pt_n(y;a,b,c|q^2) f(y) dm^\hf(y;a,b,c|q^2) \right] \ e_{n} \tensor e_{n-p} , & p \leq 0,\\
\displaystyle \sum_{n=0}^\infty \left[ \int_{-1}^1 \Pt_n(y;a,b,c|q^2) f(y) dm^\hf(y;a,b,c|q^2)\right]\ e_{n+p} \tensor e_n , & p \geq 0,
\end{cases}
\end{equation}
where $a$, $b$ and $c$ are given by \eqref{abc}.
\begin{proof}
First we concentrate on the case $k_1-k_2 \geq -\hf$ and
$k_1+k_2\geq \hf$, then $dm^\hf(\cdot;a,b,c|q^2)$ only has an
absolutely continuous part, which we denote by
$w(\cos\te;a,b,c|q^2)d\te$. We have
\[
w(\cos\te;a,b,c|q^2)=\sqrt{\frac{(q^2,ab,ac,bc;q^2)_\infty}{2\pi}}
\left| \frac{ (e^{2i\te};q^2)_\infty } { (ae^{i\te}, be^{i\te},
ce^{i\te};q^2)_\infty } \right|.
\]
Observe that
\begin{equation} \label{eq w}
w(\cos\te;aq^2,b,c|q^2) = \sqrt{ \frac{
(1-ae^{i\te})(1-ae^{-i\te}) } {(1-ab)(1-ac) }}\,
w(\cos\te;a,b,c|q^2).
\end{equation}

We define a unitary operator $\La$ by
\begin{equation} \label{La}
\La (e_{n_1} \tensor e_{n_2}) =
\begin{cases}
\displaystyle (-1)^{n_1-n_2}\int_0^\pi
\Pt_{n_1}(\cos\te;a,b,c|q^2)w(\cos\te;a,b,c|q^2) e_{n_1-n_2-L}d\te
 , & n_1 \leq n_2, \\ \displaystyle
\int_0^\pi \Pt_{n_2}(\cos\te;a,b,c|q^2)w(\cos\te;a,b,c|q^2)
e_{n_1-n_2-L}d\te
 , & n_1 \geq n_2,
\end{cases}
\end{equation}
where $e_{n-m-L}$ is an orthonormal basisvector for the
representation space of the principal unitary series, where
$e^{i\te}=q^{2i\rho}$, and $a,b,c$ are the parameters as in
\eqref{param}. We prove that $\La$ intertwines the action of $A$,
$B$, $C$, $D$ in the tensor product representation with the action
in the direct integral representation, i.e.
\begin{equation} \label{intertw}
\La \circ \pi_{k_1}^+ \tensor\pi_{k_2}^- \big( \De(Y) \big)=
\dirint \pi^P_{\rho,\eps}(Y) d \rho \circ \La, \qquad Y=A,B,C,D.
\end{equation}

We use \eqref{comult}, \eqref{pos} and \eqref{neg} to determine
the action of $B$ in the tensor product. For $n_2=n$ and
$n_1-n_2=p>0$
\begin{align*}
&(q^{-1}-q)  \La \left( \pi^+_{k_1} \tensor \pi_{k_2}^- \big(
\De(B) \big) e_{n+p} \tensor e_{n} \right) = \\ &
q^{k_2-k_1-p-\hf} \sqrt{ (1-q^{2n+2p+2})(1-q^{4k_1+2n+2p})}
\int_{0}^\pi \Pt_n(\cos\te;aq^2,b,c|q^2) w(\cos\te;aq^2,b,c|q^2) e_{p+1-L} d\te \\
&-  q^{k_1-k_2+p+\hf} \sqrt{ (1-q^{2n})(1-q^{2k_2+2n-2})}
\int_{0}^\pi \Pt_{n-1}(\cos\te;aq^2,b,c;q^2)
w(\cos\te;aq^2,b,c;q^2) e_{p+1-L} d\te
\\ &=   q^{k_2-k_1-p-\hf} \int_{0}^\pi |1-q^{2(k_1-k_2+p+\hf)}
e^{i\te}| \Pt_n(\cos\te;a,b,c|q^2) w(\cos\te;a,b,c|q^2) e_{p+1-L} d\te \\
&= (q^{-1}-q) \dirint \pi^P_{\rho,\eps} (B) d\rho \circ
\La(e_{n+p} \tensor e_n).
\end{align*}
Here we use \eqref{eq w} and the contiguous relation
\begin{equation} \label{contig2}
\begin{split}
&\sqrt{(1-ab)(1-ac)} \Pt_n(\cos\te;a,b,c|q^2) = \\ &
\sqrt{(1-abq^{2n})(1-acq^{2n})} \Pt_n(\cos\te;aq^2,b,c|q^2) -a
\sqrt{(1-q^{2n})(1-bcq^{2n-2}) } \Pt_{n-1}(\cos\te;aq^2,b,c|q^2).
\end{split}
\end{equation}
This relation can be verified by expanding
\[
\sqrt{(1-ab)(1-ac)} \Pt_n(\cos\te;a,b,c|q^2)  = \sum_{j=0}^{n} c_j
\Pt_j(\cos\te;aq^2,b,c|q^2),
\]
where
\[
c_j = \int_{\R} \sqrt{(1-ab)(1-ac)} \Pt_n(y;a,b,c|q^2)
\Pt_j(y;aq^{2},b,c|q^2) dm(y;aq^{2},b,c|q^2).
\]
From \eqref{eq w} and the orthogonality relations for $\Pt_n$
follows that $c_j=0$ for $j<n-1$. We determine $c_n$ and $c_{n-1}$
from the orthogonality relations and the leading coefficient $lc$
of $\Pt_n(\cos\te;a,b,c|q^2)$,
\[
lc = \frac{2^n}{\sqrt{(q^2,ab,ac,bc;q^2)_n }}.
\]
Relation \eqref{contig2} can also be found from theorem \ref{thm2.1} and remark \ref{rem2.2}, by substituting $a\mapsto 0$, $d \mapsto aq$ and $\al \mapsto a$.

For $p \leq 0$ we find the intertwining property of $\La$ for the
action of $B$ in the same way, using the contiguous relation
\[
\begin{split}
\frac{|1-aq^{-2} e^{i\te} |^2}{\sqrt{ (1-abq^{-2}) (1-acq^{-2}) }}
&\Pt_n(\cos\te;a,b,c|q^2) = \\ &-aq^{-2} \sqrt{
(1-q^{2n+2})(1-bcq^{2n})} \Pt_{n+1}(\cos\te;aq^{-2},b,c|q^2) \\&+
\sqrt{ (1-abq^{2n-2})(1-acq^{2n-2})}
\Pt_n(\cos\te;aq^{-2},b,c|q^2).
\end{split}
\]
This relation can be verified in the same way as \eqref{contig2}.

We can check the intertwining property of $\La$ for the action of
$C$ similarly, or we can find it using $B^*=-C$ and the fact that
$\La$ is unitary. For the action of $A$ and $D$ the intertwining
property of $\La$ follows by a straightforward computation, using
\eqref{comult}, \eqref{pos}, \eqref{neg} and \eqref{prin}.\\

For $k_1+k_2<\hf$ the measure $dm^\hf(\cdot;a,b,c|q^2)$ has  one
discrete mass point. Now the intertwining operator $\La$ is
defined similarly to \eqref{La}, where in the discrete mass point
$e_{n_1-n_2-L}$ is a standard orthonormal basisvector for the
representation space of the complementary series.\\

For $k_1-k_2 <-\hf$ the measure $dm^\hf(\cdot;a,b,c|q^2)$ has
finitely many discrete mass points. Again the intertwining
operator is defined similarly to \eqref{La}. Now in the discrete
mass points the orthonormal basisvectors are
$(-1)^{n_1-n_2}e_{k-n_1-n_2}$, where $e_{k-n_1-n_2}$ are the
standard basisvectors for the representation space of the negative
discrete series and $k$ is the summation index for the discrete
part of the measure as in \eqref{measure}.
\end{proof}

\begin{rem}
Note that the strange series representations do not appear in the decomposition of theorem \ref{decomp}. The decomposition in theorem \ref{decomp} looks similar to the decomposition of the tensor product of a positive and a negative
discrete series representation of the Lie algebra $\su(1,1)$, see \cite{GK}.
However for the quantum algebra $U_q\big(\su(1,1)\big)$ the action of the Casimir in the tensor product is bounded, contrary to the Lie algebra case, where the action of the Casimir in the tensor product is unbounded.
\end{rem}

\section{Overlap coefficients} \label{sec3}
In this section we consider the action of a self-adjoint element
$Y_sA$ in $U_q\big(\su(1,1)\big)$, where $Y_s$ is a twisted
primitive element. We determine (generalized) eigenvectors of this element in the discrete series, the principal unitary series and the complementary series representations. The Al-Salam and Chihara polynomials and the little $q$-Jacobi functions appear as overlap coefficients. Then we consider the action of $Y_sA$ in the tensor product. We find the
generalized ``uncoupled" eigenvectors for the tensor product
representation and the generalized coupled eigenvectors for the
direct integral representation. The Glebsch-Gordan coefficients
for the generalized eigenvectors turn out to be Askey-Wilson
polynomials. As a result we obtain two summation formulas for the orthogonal polynomials involved.

\subsection{Orthogonal functions and polynomials}
The Al-Salam and Chihara polynomials $s_n$ were first investigated by
Al-Salam and Chihara in \cite{ASC}. The polynomials $s_n$ form the
subclass of the Askey-Wilson polynomials with $c=d=0$;
\begin{equation} \label{ASC}
s_n(\cos\te)=s_n(\cos\te;a,b|q) = p_n(\cos\te;a,b,0,0|q)=a^{-n}(ab;q)_n
\ph{3}{2}{q^{-n},ae^{i\te},ae^{-i\te}}{ab,0}{q,q}.
\end{equation}
The orthonormal Al-Salam and Chihara polynomials $\s_n$ are orthonormal
with respect to the measure $dm(\cdot;a,b|q) = dm(\cdot;a,b,0,0|q)$ and
they satisfy the following three-term recurrence relation
\begin{equation} \label{rec ASC}
\begin{split}
2y\s_n(y) &= a_{n} \s_{n+1}(y) +q^n(a+b) \s_n(y) + a_{n-1} \s_{n-1}(y),
\\ a_n &= \sqrt{ (1-abq^{n-1})(1-q^n)}.
\end{split}
\end{equation}

In base $q^{-1}>1$ the moment problem corresponding to the Al-Salam
and Chihara polynomials is determinate for certain values of the
parameters. If we rewrite the recurrence relation \eqref{rec ASC}
for $\s_n(y;a,b|q^{-1})$, using
\[
\s_n(y;a,b|q^{-1}) = (-1)^n \left( \frac{q}{ab} \right)^{n/2} \sqrt{
\frac{ (q;q)_n } {(a^{-1}b^{-1};q)_n } } P_n(2y),
\]
we find the recurrence relation for $P_n(y)$;
\[
(1-q^{n+1})P_{n+1}(y) = (a+b-yq^n)P_n(y) - (ab-q^{n-1}) P_{n-1}(y).
\]
This is the form Askey and Ismail use in \cite[\S 3.12]{AI} to
determine the orthogonality relations. Without loss of generality
we assume $|a|\geq|b|$. The moment problem corresponding to the
polynomials $P_n$ is determinate if and only if $a \neq b$ and
$|q|\geq |b/a|$, cf. \cite[Thm. 3.2]{AI}. From
\cite[eq.(3.80)-(3.82)]{AI} we find that $\s_n(y;a,b|q^{-1})$ is
orthonormal with respect to the measure $d\mu(y;a,b|q^{-1})$
defined by
\[
\int_\R f(y) d\mu(y;a,b|q^{-1}) = \sum_{p=0}^\infty f(\mu(aq^{-p})) W_p
,
\]
where
\begin{equation} \label{Wp}
W_p=W_p(a,b|q^{-1}) =  \frac{ (1-a^{-2}q^{2p})
(a^{-2},a^{-1}b^{-1};q)_p (a^{-1}bq;q)_\infty } { (1-a^{-2}) (q,
a^{-1}bq;q)_p (a^{-2}q;q)_\infty } \left( \frac{b}{a} \right)^p q^{p^2}
.
\end{equation}
\begin{rem} \label{Jac pol}
The Al-Salam and Chihara polynomials in base $q^{-1}$ are closely related to the little $q$-Jacobi polynomials $p_n$. The orthonormal little $q$-Jacobi polynomials $\p_n$ are defined by, see \cite{AA}, \cite{GR}, \cite{KS},
\[
\begin{split}
\p_n(x;a,b;q) =& \sqrt{ \frac{ 1-abq^{2n+1} }{ 1-abq } \frac{ (aq, abq;q)_n } {(q; bq;q)_n } (aq)^{-n} }\, p_n(x;a,b;q), \\
p_n(x;a,b;q) =& \ph{2}{1}{q^{-n}, abq^{n+1} }{ aq }{q, qx}.
\end{split}
\]
The little $q$-Jacobi polynomials are orthonormal with respect to a discrete measure,
\[
\sum_{p=0}^\infty w_p(a,b;q) \p_n(q^p;a,b;q) \p_m(q^p;a,b;q) = \de_{nm}, \qquad
w_p(a,b;q) = \frac{ (bq;q)_p\,(aq;q)_\infty } { (q;q)_p \, (abq^2;q)_\infty } (aq)^p.
\]
The dual orthogonality relations for the Al-Salam and Chihara polynomials in base $q^{-1}$, are the orthogonality relations for the little $q$-Jacobi polynomials. In fact
\[
\s_n(\mu(aq^{-p});a,b|q^{-1}) \sqrt{W_p(a,b|q^{-1})} = \p_p(q^n;b/a,1/abq;q) \sqrt{w_n(b/a, 1/abq;q)}.
\]
This can be shown by transforming the $_3\varphi_2$-series in the definition \eqref{ASC} for the Al-Salam and Chihara polynomials in base $q^{-1}$  into a $_2\varphi_1$-series in base $q^{-1}$ using \cite[eq.(III.6)]{GR}. Next we write this expression in base $q$ by
\[
(\al;q^{-1})_n = ( \al^{-1};q)_n (-\al)^n q^{-\hf n(n-1)}.
\]
Now we have a terminating $_2\varphi_1$-series in base $q$. Finally reversing the summation and using one of Heine's transformations \cite[eq.(III.2)]{GR}, we find the desired expression.\\

\end{rem}

The little $q$-Jacobi functions $\phi_n(\mu(y);c,d,r|q)$, see
Kakehi \cite{Ka}, Koelink and Stokman \cite[\S A.2]{KSt2},
\cite{KSt1}, are defined by
\begin{equation} \label{def Jac}
\phi_n\big( \mu(y)\big)=\phi_n\big(\mu(y);c,d,r|q\big) =e^{i\psi_n}
|d|^{-n}
\sqrt{\frac{(\overline{c}q^{1+n}/\overline{d}^2\overline{r};q)_\infty }
{(q^{1+n}/r;q)_\infty }} \ph{2}{1}{dy, d/y}{c}{q,rq^{-n}}, \quad n \in
\Z,
\end{equation}
where $\psi_n \in \R$ is such that
\[
\psi_{n+1}-\psi_n = \arg\left(\overline{d}(1-\frac{q^{n+1}}{r})\right)
= \arg\left( \overline{d}
(1-\frac{\overline{c}q^{n+1}}{\overline{d}^2\overline{r}})\right).
\]
Note that we use a slightly different definition than in \cite{KSt1},
\cite{KSt2}. The little $q$-Jacobi functions satisfy the recurrence
relation
\begin{equation} \label{rec phi}
\begin{split}
 2x \phi_n(x) &= a_{n} \phi_{n+1}(x) +
\frac{q^n(c+q)}{dr} \phi_n(x) + a_{n-1} \phi_{n-1}(x),\\ a_n &= \sqrt{
(1-\frac{q^{n+1}}{r})(1-\frac{cq^{n+1}}{d^2r}) }.
\end{split}
\end{equation}
For $0<c \leq q^2$, $|c/d|<1$, $|d|<1$, $r, c/d^2r \notin q^\Z$ and (1)
$\overline{r}c=d^2r$ or (2) $r>0$, $c \neq d^2$ and $rq^{k_0+1}
< c/d^2
< rq^{k_0}$, where $k_0 \in \Z$ is such that $1<rq^{k_0}<q^{-1}$, the
little $q$-Jacobi functions $\phi_n(x)$ are orthonormal with respect to
the measure $d\nu(\cdot;c,d,r|q)$ given by
\begin{equation} \label{meas2}
\begin{split}
\int_\R & f(x) d\nu(x;c,d,r|q) = \\& \frac{1}{2\pi}\int_{0}^\pi
f(\cos\te) \left| \frac{(e^{2i\te},c,r,q/r;q)_\infty} { (ce^{i\te}/d,
de^{i\te}, dre^{i\te}, qe^{i\te}/dr;q)_\infty } \right|^2 d\te \\& +
\sum_{\substack{j \in \Z\\ |q^{-j}/dr|>1}} f\big(\mu(q^{-j}/dr)\big)
 \frac{|(c,r,q/r;q)_\infty|^2}{(q;q)_\infty^2} \frac{
(1-d^2r^2q^{2j}) (dr)^{-2j-2} q^{-j(j+1)} } {
(crq^{j},cq^{-j}/d^2r,d^2rq^{j},q^{-j}/r;q)_\infty } ,
\end{split}
\end{equation}
Note that more general conditions for orthogonality can be given, see
\cite{KSt2}. The conditions given here are related to the action of a
twisted primitive element in the principal unitary series (1) and in
the complementary series (2). This will be explained in the next
subsection.

\subsection{Eigenvectors for $\boldsymbol{\pi^+}$, $\boldsymbol{\pi^-}$,
 $\boldsymbol{\pi^P}$ and $\boldsymbol{\pi^C}$}
For $s \in \R \setminus\{0\}$ we define an element
$Y_s$ by
\[
Y_s = q^{\frac{1}{2}}B-q^{-\frac{1}{2}}C +
\frac{s+s^{-1}}{q^{-1}-q}(A-D) \in U_q\big(\su(1,1)\big).
\]
$Y_s$ is a twisted primitive element, i.e. $\De(Y_s) = A \tensor Y_s + Y_s \tensor D$. Twisted primitive elements are elements which are much like Lie algebra elements, see e.g. \cite{Koo1}. We consider the action of
\[
Y_sA = q^{\frac{1}{2}}BA-q^{-\frac{1}{2}}CA +
\frac{s+s^{-1}}{q^{-1}-q}(A^2-1),
\]
which is a self-adjoint element in $U_q\big(\su(1,1)\big)$. We start
with the action in the discrete series representations.
\begin{prop} \label{prop1}
The operator $\Te^+$ defined by
\begin{align*}
\Te^+ : \lt(\Z_{\geq0}) &\rightarrow L^2(\R, dm(\cdot;sq^{2k},
q^{2k}/s|q^2)\\ e_n &\mapsto \s_n(\cdot;sq^{2k}, q^{2k}/s|q^2)
\end{align*}
is unitary and intertwines $\pi^+_k(Y_sA)$ acting on $\lt(\Z_{\geq0})$
with $(q^{-1}-q)^{-1}M_{2x-2\mu(s)}$ acting on \\ $L^2\big(\R,
dm(x;sq^{2k}, q^{2k}/s|q^2)\big)$.
\end{prop}
This is proposition $4.1$ in Koelink and Van der Jeugt \cite{KJ}. It is
proved by showing that $\pi^+_k(Y_sA)$ is a bounded self-adjoint
Jacobi operator, corresponding to the recurrence relation \eqref{rec
ASC} for the Al-Salam and Chihara polynomials. Note that the spectrum of
$\pi^+_k(Y_sA)$ consists of a bounded continuous part and a (possibly empty)
finite discrete part. The proof of the next
proposition is completely analogous, using an unbounded, essentially
self-adjoint Jacobi operator.
\begin{prop} \label{prop2}
For $|s|\geq q^{-1}$ the operator $\Te^-$ defined by
\begin{align*}
\Te^- : \lt(\Z_{\geq0}) &\rightarrow L^2\big(\R,d\mu(\cdot;
sq^{-2k},q^{-2k}/s|q^{-2}) \big)\\ e_n &\mapsto (-1)^n
\s_n(\cdot;sq^{-2k}, q^{-2k}/s|q^{-2})
\end{align*}
is unitary and intertwines $\pi^-_k(Y_sA)$ acting on $\lt(\Z_{\geq0})$
with $(q^{-1}-q)^{-1}M_{2\mu(q^{2k+2p}/s)-2\mu(s)}$ acting on
$L^2\big(\R,d\mu(\cdot; sq^{-2k},q^{-2k}/s|q^{-2}) \big)$.
\end{prop}
Observe that the spectrum of $\pi^-_k(Y_sA)$ is discrete and unbounded.
From here on we assume $|s| \geq q^{-1}$ in order to make the
moment problem corresponding to the Al-Salam and Chihara
polynomials in base $q^{-2}$ a determinate moment problem.

From proposition \ref{prop1} we conclude that for $x \in [-1,1]$
\begin{equation} \label{v+}
v^+(x) = \sum_{n=0}^\infty \s_n(x;sq^{2k}, q^{2k_1}/s|q^2) e_n
\end{equation}
is a generalized eigenvector of $\pi^+_k(Y_sA)$ for the eigenvalue
$2(x-\mu(s))/(q^{-1}-q)$. If $x$ is in the discrete part of the support
of the measure $dm$ of proposition \ref{prop1}, then $v^+(x)$ is an eigenvector
of $\pi^+_k(Y_sA)$. Further we conclude from proposition \ref{prop2}
that
\begin{equation} \label{v-}
v^-(p) = \sum_{n=0}^\infty (-1)^n
\s_n\big(\mu(q^{2k+2p}/s);sq^{-2k},q^{-2k}/s|q^{-2}\big) e_n
\end{equation}
is an eigenvector of $\pi^-_k(Y_sA)$ for the eigenvalue
$2(\mu(q^{2k+2p}/s)-\mu(s))/(q^{-1}-q)$.

Next we consider the action of $Y_sA$ in the principal unitary series
and in the complementary series.
\begin{prop} \label{prop3}
The operator $\Te^P$ defined by
\begin{align*}
\Te^P : \lt(\Z) &\rightarrow L^2(\R, d\nu(\cdot; q^2/s^2,
q^{2i\rho+1}/s, q^{1-2\eps-2i\rho}|q^2) \\ e_n &\mapsto
\phi_n(\cdot;q^2/s^2, q^{2i\rho+1}/s, q^{1-2\eps-2i\rho}|q^2)
\end{align*}
is unitary and intertwines $\pi^P_{\rho,\eps}(Y_sA)$ acting on
$\lt(\Z)$ with $(q^{-1}-q)^{-1}M_{2x-2\mu(s)}$ acting on \\ $L^2\big(\R,
d\nu(x; q^2/s^2, q^{2i\rho+1}/s, q^{1-2\eps-2i\rho}|q^2)\big)$.
\end{prop}
Note that the spectrum of $\pi^P_{\rho,\eps}(Y_sA)$ consists of a bounded continuous part and an unbounded discrete part. Proposition \ref{prop3} generalizes the discussion in \cite[\S6]{KSt2}.
\begin{proof}
Using \eqref{prin} we see that $\pi^P_{\rho,\eps}(Y_sA)$ is a doubly
infinite Jacobi operator, see \cite{MR}, \cite{Koe}, corresponding to
the recurrence relation for the little $q$-Jacobi functions \eqref{rec
phi}. The result then follows from the spectral decomposition, which is
equivalent to the orthonormality relations for the little $q$-Jacobi
functions.
\end{proof}
\begin{rem} \label{rem1}
(i) Observe that the support of the measure $d\nu$ in proposition \ref{prop3}
does not depend on $\rho$, so that the spectral value $2\big(x-\mu(s)\big)/(q^{-1}-q)$ also does not depend on $\rho$.

(ii) For this set of parameters, a convenient expression for the
phasefactor $e^{i\psi_n}$ in the definition of the little $q$-Jacobi
functions, is
\[
e^{i\psi_n} = (-\sgn s)^n \,  \frac{(q^{1-2n-2\eps-2i\rho};q^2)_\infty
} {|(q^{1-2n-2\eps-2i\rho};q^2)_\infty | }.
\]
If we now transform the $_2\varphi_1$-series in the definition of
the little $q$-Jacobi function into a $_2\varphi_2$-series, using
\cite[eq.(III.4)]{GR}, we find
\[
\begin{split}
\phi_n(\mu(y);q^2/s^2,& q^{2i\rho+1}/s, q^{1-2\eps-2i\rho}|q^2) = \\& \left(-\frac{s}{q}\right)^n \frac{ (yq^{2-2\eps-2n}/s;q^2)_\infty }{ | (q^{1-2\eps-2n-2i\rho};q)_\infty |} \ph{2}{2}{yq^{1+2i\rho}/s, yq^{1-2i\rho}/s} {q^2/s^2,yq^{2-2\eps-2n}/s}{q^2, q^{2-2\eps-2n}/ys}.
\end{split}
\]
So we see that $\phi_n(x;q^2/s^2,
q^{2i\rho+1}/s, q^{1-2\eps-2i\rho}|q^2)$ is symmetric in
$q^{2i\rho}$ and $q^{-2i\rho}$. To stress that $\phi_n(x;q^2/s^2,
q^{2i\rho+1}/s, q^{1-2\eps-2i\rho}|q^2)$, besides a function of $x$, is also a function of $\rho$, we will also use the notation $\phi_n(x;\cos\psi)$, where $e^{i\psi}=q^{2i\rho}$.

(iii) Formally, for $q^{2i\rho} = q^{2k_1-2k_2-2j+1}>1$, i.e.~$y=\mu(q^{2i\rho})$ is a discrete mass point of the measure $dm$ for the continuous dual $q$-Hahn polynomials corresponding to a negative discrete series reprensentation, the little $q$-Jacobi function becomes an Al-Salam and Chihara polynomial in base $q^{-2}$. Indeed, for these values of $\rho$ the $_2\varphi_1$-series in the definition of the little $q$-Jacobi functions \eqref{def Jac} can be transformed into a terminating $_2\varphi_1$-series.
\end{rem}

Proposition \ref{prop3} states that for $x \in [-1,1]$
\begin{equation} \label{vP}
v^P(x) = \sum_{n=-\infty}^\infty \phi_n(x;q^2/s^2, q^{2i\rho+1}/s,
q^{1-2\eps-2i\rho}|q^2) e_n
\end{equation}
is a generalized eigenvector of $\pi^P_{\rho,\eps}(Y_sA)$ for the
eigenvalue $2(x-\mu(s))/(q^{-1}-q)$. If $x$ is in the discrete part of the support of the measure $d\nu$ of proposition \ref{prop3}, then $v^P(x)$ is a real eigenvector of $\pi^P_{\rho,\eps}(Y_sA)$. Also proposition \ref{prop3}
shows that $\pi^P_{\rho,\eps}(Y_sA)$ is essentially self-adjoint for $|s| \geq q^{-1}$.

The result for the complementary series representations follows
directly from proposition \ref{prop3} by replacing $-\hf+i\rho$ by
$\la$.

\subsection{Coupled and uncoupled eigenvectors}
Next we consider the action of $Y_sA$ in the tensor product
representation $\pi^+_{k_1} \tensor \pi^-_{k_2}$. Denote
\begin{equation} \label{eq:F}
F_{n_1,n_2}(x,p) =
(-1)^{n_2}\s_{n_1}(x;sq^{2k_1-2k_2-2p},q^{2k_1+2k_2+2p}/s|q^2) \s_{n_2}(\mu(q^{2k_2+2p}/s);
q^{-2k_2}s, q^{-2k_2}/s |q^{-2}).
\end{equation}
and let $dm(x,p)$ denote the measure
\begin{equation} \label{dm(xp)}
dm(x;sq^{2k_1-2k_2-2p},q^{2k_1+2k_2+2p}/s|q^2)d\mu(x_2;q^{-2k_2}s,
q^{-2k_2}/s |q^{-2}).
\end{equation}
Recall that $d\mu(x_2)$ is a discrete measure, with mass points in
$x_2=\mu(q^{2k_2+2p}/s)$, $p \in \Z_{\geq 0}$. Observe that
the number of discrete points of the
orthogonality measure $dm(x)$ for the Al-Salam and Chihara polynomials $\s_{n_1}$
depends on $p$. For $p \rightarrow \infty$ the number of discrete mass points of $dm(x)$ tends to infinity.
\begin{prop} \label{prop4}
The operator $\Ups$ defined by
\begin{align*}
\Ups : \lt(\Z_{\geq 0}) \tensor \lt(\Z_{\geq 0}) &\rightarrow
L^2\left(\R \times \Z_{\geq0}, dm(x,p)\right) \\ e_{n_1}
\tensor e_{n_2} &\mapsto F_{n_1,n_2}(x,p)
\end{align*}
is unitary and intertwines $\pi^+_{k_1} \tensor \pi^-_{k_2} \big(
\De(Y_sA) \big)$ acting on $\lt(\Z_{\geq 0}) \tensor \lt(\Z_{\geq 0})$
with $(q^{-1}-q)^{-1}M_{2x-2\mu(s)}$ acting on $L^2\big(\R \times \Z_{\geq 0}, dm(x,p)\big)$.
\end{prop}
\begin{proof}
We have $\De(Y_sA)= A^2 \tensor Y_sA + Y_sA \tensor 1$, so the
action of $Y_sA$ is
\[
\begin{split}
(q^{-1}-q) \pi^+_{k_1} \tensor \pi^-_{k_2} \big( \De(Y_sA) \big)
e_{n_1} \tensor e_{n_2} &=\\ (q^{-1}-q) \left[\pi^+_{k_1} (A^2)
e_{n_1} \right] \tensor & \left[ \pi^-_{k_2}(Y_sA) e_{n_2}\right]
+ (q^{-1}-q) \left[\pi^+_{k_1} (Y_sA) e_{n_1}\right] \tensor
e_{n_2}
\end{split}
\]
Define for fixed $p$ an operator $\Ups_{p}$ by
\[
\Ups_{p} (e_{n_1} \tensor e_{n_2}) =
(-1)^{n_2}\s_{n_2}(\mu(q^{2k_2+2p}/s); q^{-2k_2}s, q^{-2k_2}/s |q^{-2})
\, e_{n_1},
\]
then from \eqref{pos} and proposition \ref{prop2} we find that
$\Ups_{p}\Big[(q^{-1}-q) \pi^+_{k_1} \tensor \pi^-_{k_2} \big(
\De(Y_sA) \big)+2\mu(s)\Big]$ acts as a Jacobi operator that can be
identified with the three-term recurrence relation for the
Al-Salam and Chihara polynomials \eqref{rec ASC}. Then the
proposition follows in the same way as proposition \ref{prop1}.
\end{proof}
Proposition \ref{prop4} states that for $x \in [-1,1]$
\begin{equation} \label{v(xp)}
v(x,p) =  \sum_{n_1=0}^\infty
\sum_{n_2=0}^\infty F_{n_1,n_2}(x,p)\ e_{n_1} \tensor e_{n_2}
\end{equation}
is a generalized eigenvector of $\pi^+_{k_1} \tensor \pi^-_{k_2} \big(
\De(Y_sA) \big)$ for eigenvalue $2\big(x_{(p)}-\mu(s)\big)/(q^{-1}-q)$.
If $x$ is in the discrete part of the measure
$dm(x;sq^{2k_1-2k_2-2p},q^{2k_1+2k_2+2p}/s|q^2)$, then $v(x,p)$ is an
eigenvector of $\pi^+_{k_1} \tensor \pi^-_{k_2} \big(
\De(Y_sA) \big)$.
We call $v(x,p)$ the uncoupled
eigenvector of $\pi^+_{k_1} \tensor \pi^-_{k_2} \big(\De(Y_sA) \big)$,
even though
\[
v(x,p) \neq v^+(x) \tensor v^-(p),
\]
due to the noncocommutativaty of the coproduct $\De$. Proposition \ref{prop4} also implies that $\pi^+_{k_1}\tensor\pi^-_{k_2} \big( \De(Y_sA) \big)$ is essentially self-adjoint.\\

Next we want to determine how $\Ups$ acts on elements in the
representation space of the direct integral representation $\int^\oplus \pi^P_{\rho,\eps} d\rho$. If there do not appear discrete terms in the decomposition of theorem \eqref{decomp}, we define for appropriate functions $g$ an operator $\Ups_g$ by
\[
\begin{split}
\Ups_g : L^2(-1,1) \tensor \lt(\Z) &\cong
\sideset{}{^\oplus}\int\limits_{0\ }^{\ \pi}\lt(\Z) d\psi
\rightarrow \sideset{}{^\oplus}\int\limits_{0\ }^{\ \pi}
L^2\big(\R, d\nu(\cdot;q^2/s^2, qe^{i\psi}/s,
q^{1-2\eps}e^{-i\psi}|q^2)\big) d\psi
\\ f \tensor e_n &\mapsto \int_0^\pi f(\cos\psi) g(\cos\psi)
\phi_n(x;q^2/s^2, qe^{i\psi}/s, q^{1-2\eps}e^{-i\psi}|q^2) d\psi,
\end{split}
\]
where we put $e^{i\psi}= q^{2i\rho}$. If discrete terms appear in the decomposition in theorem \ref{decomp}, we must add discrete terms in the definition of $\Ups_g$. We leave this to the reader. We use proposition \ref{prop3} to see that $\Ups_g$ intertwines $\int^\oplus \pi^P_{\rho,\eps}(Y_sA) d\rho$ with $(q^{-1}-q)^{-1}M_{2x-2\mu(s)}$ for any function $g$. Now from the Clebsch-Gordan decomposition \eqref{CG decomp} and proposition \ref{prop4} we find that there exists a unique function $g$, such that $\Ups = \Ups_g$, or equivalently, for $n= \min\{n_1,n_2\}$,
\begin{equation} \label{int g}
F_{n_1,n_2}(x,p) =  \int_\R
g(y)\Pt_{n}(y) \phi_{n_1-n_2-L}(x;y)
dm^\hf(y),
\end{equation}
where $dm$ is the orthogonality measure for the continuous dual $q$-Hahn polynomials $\Pt_n$ with parameters as in theorem \ref{decomp}. We write $\phi_n(x;y)$ to stress the fact that the little $q$-Jacobi function depends on $x$ as well as $y=\cos\psi=\mu(q^{2i\rho})$, see remark \ref{rem1}(ii). Recall that $\eps= {k_1-k_2+L}$. Observe that $g$ is the Clebsch-Gordan coefficient for the bases of (generalized) eigenvectors $v(x,p)$ and $v^P(x)$. Therefore $g$ does not depend on $n_1$ and $n_2$. Also $g$ is uniquely determined by the the choice of bases and the intertwiner in theorem \ref{decomp}.

\begin{figure}[t]
\includegraphics[width=15cm, height=6cm]{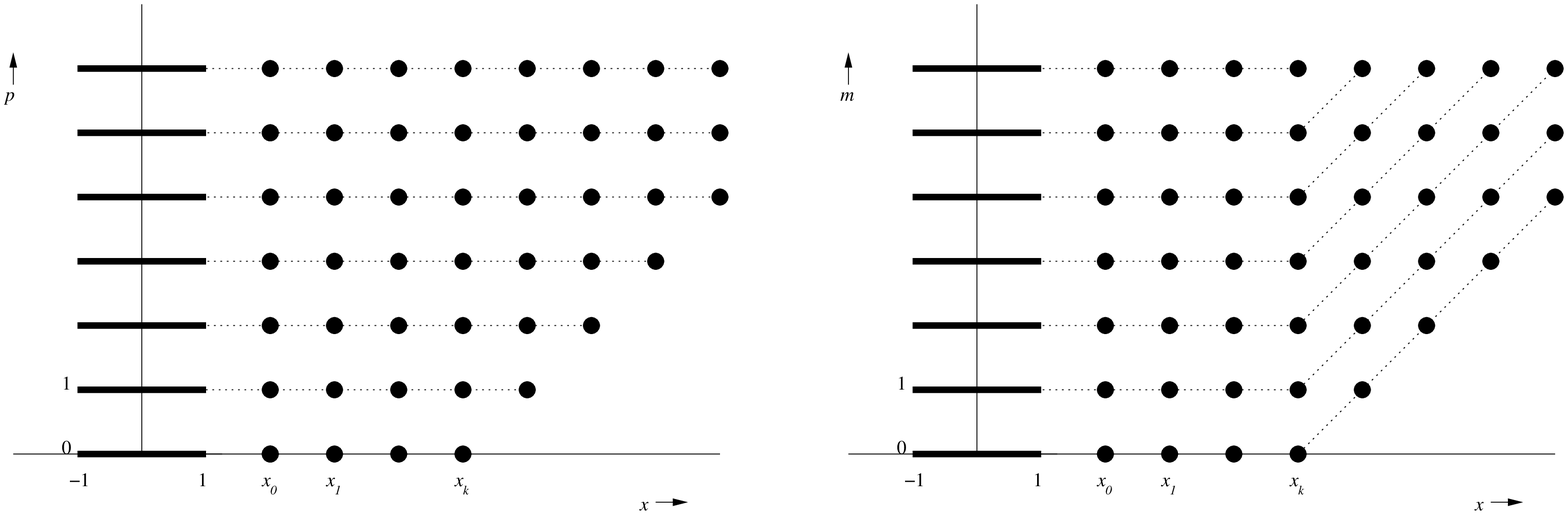}
\caption{Left: a simplified picture of the support of $dm(x,p)$. Right: the same picture with variables $x$ and $m$.} \label{fig1}
\end{figure}

The determination of the Clebsch-Gordan coefficient $g$ is the crucial step in this paper. To sketch the idea we have depicted the support $\mathcal D$ of the measure $dm(x,p)$ of proposition \ref{prop4} in figure \ref{fig1} in two ways. On the left we have drawn horizontal lines to depict the orthogonality of $F_{n_1,n_2}$ on $\mathcal D$, cf. \eqref{eq:F}, \eqref{dm(xp)}. On the right we have split up $\mathcal D$ in broken lines, and each broken line is the support of the measure $d\nu$ of proposition \ref{prop3}. So the little $q$-Jacobi functions form an orthogonal set of functions on each broken line. It remains to find the orthogonal functions in the $m$-direction to complement the little $q$-Jacobi functions to an orthogonal basis on $\mathcal D$. This gives the function $g$. Because of all the discrete mass points a lot of bookkeeping is necessary. \\

Let $|x_0|<|x_1|<|x_2|<\ldots$ denote the discrete mass points of the orthogonality measure for the little $q$-Jacobi polynomials $d\nu(x; q^2/s^2, qe^{i\psi}/s, q^{1-2\eps}e^{-i\psi}|q^2)$. So $x_0 = \mu(sq^{2k_1-2k_2-2j})$, $\dots$ , $x_l = \mu(sq^{2k_1-2k_2-2j-2l})$, where $j$ is the smallest integer such that $|s|q^{2k_1-2k_2-2j}>1$. Further let $k+1$ be the number of discrete mass points of the orthogonality measure $dm(\cdot;sq^{2k_1-2k_2}, q^{2k_1+2k_2}/s|q^2)$ for the Al-Salam and Chihara polynomials. Observe that these discrete points are exactly the points $x_0, \ldots, x_k$, so $j=-k$. We define a measure $dM(\cdot;x,p)$ by
\[
dM(y;x,p)dm(x,p) = dm(y;a,b,c,d|q^2) d\nu(x; q^2/s^2, qe^{i\psi}/s, q^{1-2\eps}e^{-i\psi}|q^2), \quad y=\cos\psi,
\]
where $dm(x,p)$ is given by \eqref{dm(xp)}. For $x=\mu(u)\in [-1,1]$ and for $x=\mu(u)=x_l$, $l\leq k$ the parameters $a$, $b$, $c$ and $d$ are given by
\begin{equation} \label{abcd1}
a=qu/s, \quad  b= q/us, \quad c= q^{2k_2-2k_1+1}, \quad d=
q^{2k_1+2k_2-1},
\end{equation}
and for $x=x_l$, $l >k$,
\begin{equation} \label{abcd2}
a=q^{2k_1-2k_2+1}, \quad  b= q^{2k_2-2k_1+2l-2k+1}/s^2, \quad c= q^{2k_2-2k_1+2l-2k+1}, \quad d=q^{2k_1+2k_2-1}.
\end{equation}
We denote the value of the measure $dM(\cdot;x,p)$ in a points y by $W(y;x,p)$. Finally let $m \in \Z_{\geq 0}$ be the number
\begin{equation} \label{deg m}
m=
\begin{cases}
p, & \text{if }x \in [-1,1] \textrm{ or } x=x_l,\ l \leq k, \\
p+k-l, & \text{if }x=x_l, \ l > k.
\end{cases}
\end{equation}
This relation is depicted in figure \ref{fig1}. Now we claim that the function $g$ is given by
\[
g(y) = e^{i\al} \p_m(y;a,b,c,d|q^2) \sqrt{W(y;x,p)},
\]
where $\al \in \R$ is fixed. To stress that the Askey-Wilson polynomials with parameters \eqref{abcd1}, \eqref{abcd2} depend on $x=\mu(u)$, we use the notation
\[
\p_m(y;a,b,c,d|q^2)= \p_m(y;x).
\]

To verify our claim we prove the following proposition.
\begin{prop} \label{prop5}
Let $m \in \Z_{\geq 0}$ be as in \eqref{deg m}, then for $n= \min\{n_1,n_2\}$
\begin{equation} \label{int1}
F_{n_1,n_2}(x,p)=
(-\sgn s)^L \int_\R \p_m(y;x) \Pt_{n}(y)
\phi_{n_1-n_2-L}(x;y) dm^\hf(y) dM^\hf(y;x,p).
\end{equation}
\end{prop}
\begin{proof}
First observe that it is enough to prove the proposition for $n_1=n_2=0$, then \eqref{int1} follows directly from \eqref{int g}.

We use the generating function for the Askey-Wilson polynomials from theorem \ref{thm2.1}. We write the $_8\varphi_7$-series as a sum of two $_4\varphi_3$-series, using \cite[eq.(III.36)]{GR}. Now first we let $t \rightarrow 0$ and then let $r \rightarrow 0$, then the second term vanishes and  the first term becomes a $_2\varphi_2$-series. So we find after using Jackson's transformation \cite[eq.(III.4)]{GR}
\begin{equation} \label{expan}
\frac{ (abcd, ce^{-i\psi};q)_\infty } { (ac, bc;q)_\infty } \ph{2}{1}
{ae^{i\psi} , be^{i\psi} }{ ab } {q, ce^{-i\psi}} = \sum_{n=0}^\infty
\frac{ (abcd;q)_{2n} (abc)^n q^{n(n-1)} p_n(\cos \psi ;a,b,c,d|q) } {
(q,ab,ac,bc,abcdq^{n-1};q)_n } .
\end{equation}
We leave it to the reader to verify that it is allowed to interchange limits and summations. Formula \eqref{expan} can also be proved directly using \eqref{AW pol} and a limiting case of Jackson's summation \cite[eq.(II.20)]{GR}. We choose the parameters $a$, $b$, $c$, $d$ as in \eqref{abcd1}, then \eqref{expan} gives an expansion of the little $q$-Jacobi function $\phi_{-L}(\mu(u);q^2/s^2,qe^{i\psi}/s,
q^{2k_2-2k_1-2L+1}e^{-i\psi}|q^2)$ in terms of Askey-Wilson polynomials of argument $\cos\psi$.

We consider first the case $x = \mu(u) \in [-1,1]$ and $k_1+k_2\geq \hf$. We use \eqref{on AW} to write the Askey-Wilson polynomials in the orthonormal form. We multiply for fixed $p \in \Z_{\geq 0}$ both sides of equation \eqref{expan} with $\p_p(\cos\psi;a,b,c,d|q^2)$, where $a$, $b$, $c$, $d$ are given by \eqref{abcd1}, and integrate against the orthogonality measure $dm(\cos\psi;a,b,c,d|q^2)$, see \eqref{measure}. Using the orthogonality of the Askey-Wilson polynomials, we obtain
\begin{equation} \label{eqA}
\begin{split}
\frac{1}{2\pi} \int_0^\pi & \ph{2}{1}{que^{i\psi}/s, qe^{i\psi}/us
}{q^2/s^2}{q^2, q^{2k_2-2k_1+1}e^{-i\psi}} \p_p(\cos\psi;qu/s, q/us,
q^{2k_2-2k_1+1}, q^{2k_1+2k_2-1}|q^2) \\&\times  \frac{ (q^2, q^2/s^2,
uq^{2k_2-2k_1+2}/s, q^{2k_2-2k_1+2}/us, uq^{2k_1+2k_2}/s,
q^{2k_1+2k_2}/us, q^{4k_2};q^2)_\infty }{ (q^{4k_2+2}/s^2;q^2)_\infty }
\\& \times \left| \frac{ (e^{2i\psi};q^2)_\infty } { (que^{i\psi}/s,
qe^{i\psi}/us, q^{2k_2-2k_1+1}e^{i\psi},
q^{2k_1+2k_2-1}e^{i\psi};q^2)_\infty } \right|^2 d\psi \\ =& \sqrt{
\frac{ (q^{4k_2+2}/s^2;q^2)_{2p} (uq^{2k_1+2k_2}/s, q^{2k_1+2k_2}/us,
q^{4k_2};q^2)_p } { (q^2, q^2/s^2, uq^{2k_2-2k_1+2}/s,
q^{2k_2-2k_1+2}/us, q^{4k_2+2p}/s^2;q^2)_p }} (q^{2k_2-2k_1+1}/s^2)^p
q^{2p^2}.
\end{split}
\end{equation}
By a straightforward calculation we see that \eqref{eqA} is the same as \eqref{int1} with $n_1=n_2=0$ and $m=p$.

For $x = \mu(u) \in [-1,1]$ and $k_1+k_2 < \hf$ the calculations proceed analogously. In this case $dm(\cdot;a,b,c,d|q)$, where $a$, $b$, $c$, $d$ are given by \eqref{abcd1}, has one discrete mass point. This discrete point corresponds to the complementary series representation occurring in the tensor product decomposition in theorem \ref{decomp}.

Next let us consider the case $x=\mu(u)=x_l$, $l \in \Z_{\geq 0}$, where $|x_0|<|x_1|<|x_2|<\ldots$ denote the discrete mass points of the measure $d\nu(\cdot; q^2/s^2, qe^{i\psi}/s, q^{2k_2-2k_1-2L+1}e^{-i\psi}|q^2)$. Note that the sum for the discrete part of the measure $d\nu$, see \eqref{meas2}, starts at index $-(k+1)$. In this case we find from \eqref{abcd1},
\[
a= q^{2k_1-2k_2+2k-2l+1}, \qquad b=q^{2k_2-2k_1+2l-2k+1}/s^2.
\]
Note that $|b|<1$ for all $l$. So the orthogonality measure $dm(\cdot;a,b,c,d|q^2)$ for the Askey-Wilson polynomials only has discrete mass points for $a>1$ and $d>1$.

For $l \leq k$ we can just repeat the proof as for $x \in [-1,1]$. So we now assume $l > k$. First consider the case $m=0$. From figure \ref{fig1} we see that $x_l$ is not inside the support of $dm(\cdot;sq^{2k_1-2k_2}, q^{2k_1+2k_2}/s|q^2)$, i.e. $x_l$ does not lie on the level $p=0$. So in this case $p \neq m$. From the figure it is clear that we must have $m=l-k$. For arbitrary $m$ we find in the same way $m=p+k-l$. We now prove that this is indeed true.

Again we repeat the proof as for $x \in [-1,1]$, but only for $p \geq l-k$. We still use here the parameters $a$, $b$, $c$, $d$ given by \eqref{abcd1}, then we find an equation similar to \eqref{eqA}. We put $n=k-l$ and we use \eqref{AW pol} to rewrite the $_4\varphi_3$-series in the Askey-Wilson polynomial $p_p(\cos\psi;a,b,c,d|q)$;
\begin{equation} \label{eq:shift}
\begin{split}
&\frac{(q^{2-2n};q^2)_\infty }{ (q^{2-2n};q^2)_p} p_p(\cos\psi;q^{2k_1-2k_2-2n+1}, q^{2k_2-2k_1+2n+1}/s^2, q^{2k_2-2k_1+1}, q^{2k_1+2k_2-1}|q^2)  \\
= &(q^{2k_1-2k_2-2n+1})^{-p} (q^2/s^2, q^{4k_1-2n};q^2)_p\\
&\times \sum_{j=n}^p \frac{ (q^{-2p}, q^{4k_2+2p}/s^2, q^{2k_1-2k_2-2n+1}e^{i\psi}, q^{2k_1-2k_2-2n+1}e^{-i\psi};q^2)_j\, q^{2j} } { (q^2, q^2/s^2, q^{4k_1-2n} ;q^2)_j } (q^{2-2n+2j};q^2)_\infty \\
= &(q^{2k_1-2k_2-2n+1})^{-p} (q^2/s^2, q^{4k_1-2n};q^2)_p\\
& \times \sum_{i=0}^{p-n} \frac{(q^{-2p}, q^{4k_2+2p}/s^2, q^{2k_1-2k_2-2n+1}e^{i\psi}, q^{2k_1-2k_2-2n+1}e^{-i\psi};q^2)_{n+i}\, q^{2n+2i} } { (q^2, q^2/s^2, q^{4k_1-2n} ;q^2)_{n+i} } (q^{2+2i};q^2)_\infty
\\ = &\frac{ (q^{2+2n};q^2)_\infty}{ (q^{2+2n};q^2)_{p-n} } ( q^{-2p}, q^{4k_1+2p}/s^2;q^2)_n \,| (q^{2k_1-2k_2-2n+1}e^{i\psi};q^2)_n|^2 q^{-n(2k_1-2k_2-2p-1)}\\
& \times p_{p-n}(\cos\psi; q^{2k_1-2k_2+1}, q^{2k_2-2k_1+2n+1}/s^2, q^{2k_2-2k_1+2n+1}, q^{2k_1+2k_2-1};q^2)
\end{split}
\end{equation}
So now we find an Askey-Wilson polynomial $p_{p-n}=p_{p+k-l}$ with parameters given by \eqref{abcd2}. Again a straightforward calculation gives \eqref{int1} with $n_1=n_2=0$.
\end{proof}

\begin{rem} \label{rem2}
(i) Note that discrete terms in \eqref{int1} corresponding to the sum of negative discrete series in the decomposition (see theorem \ref{decomp}) occur only when $x=x_l$. Since (at least formally) the little $q$-Jacobi polynomials become Al-Salam and Chihara polynomials in these discrete terms, this corresponds to the fact that the Al-Salam and Chihara polynomials in base $q^{-2}$  are orthogonal with respect to a discrete measure.

(ii) An interesting special case is when $n_2=p=0$. Then proposition \ref{prop5} gives a integral representation for the Al-Salam and Chihara polynomial in base $q^2$.

(iii) Let $\mathcal D \subset \R^2$ denote the support of the measure $dm(x,p)$. Then it is clear that the polynomials $F_{n_1,n_2}(x,p)$ form an orthonormal basis on $\mathcal D$. Since the functions $G_{y,n}(x,m)=p_m(y;x)\phi_{n}(x;y)$ also form an orthonormal basis on $\mathcal D$, proposition \ref{prop5} gives a connection formula between two orthonormal bases on $\mathcal D$ and the continuous dual $q$-Hahn polynomials have an interpretation as connection coefficients. Figure \ref{fig1} gives a picture of $\mathcal D$ for both orthonormal systems.

(iv) Proposition \ref{prop5} gives a formal relation between the generalized eigenvectors \eqref{v(xp)} and \eqref{vP};
\begin{equation} \label{CG ev}
v(x,p) = (- \sgn s)^L \int_{\R} \p_m(y;x) v^P(x) dM^\hf(y;x,p).
\end{equation}
This shows that the Askey-Wilson polynomials have an interpretation as Glebsch-Gordan coefficients for $U_q\big(\su(1,1)\big)$.
\end{rem}

\subsection{Summation formulas}
From proposition \ref{prop5} we can derive two summation formulas for the Al-Salam and Chihara polynomials; the first using the orthogonality of the continuous dual $q$-Hahn polynomials, the second using the orthogonality of the Askey-Wilson polynomials. For simplicity we assume $x=\cos\te \in [-1,1]$ and $y=\cos\psi \in [-1,1]$. However it is clear from the proofs how to extend the results to the general case.

\begin{thm} \label{thm1}
Let $m \in \Z$, $p \in \Z_{\geq 0}$ and $s \in \R$, with $|s| \geq q^{-1}$. For the Askey-Wilson polynomials $p_p$, the
continuous dual $q$-Hahn polynomials $P_n$ and the Al-Salam and
Chihara polynomials $s_n$, the following summation formula holds
\begin{equation} \label{sum1}
\begin{split}
\sum_{n=0}^\infty c_n s_{n+m}&(\cos\te;sq^{2k_1-2k_2-2p},
q^{2k_1+2k+2p}/s|q^2) s_{n}(\mu(q^{2k_2+2p}/s);sq^{-2k_2},q^{-2k_2}/s
|q^{-2}) \\ &  \times P_n(\cos\psi; q^{2k_1-2k_2+2m+1},q^{2k_1+2k_2-1},
q^{2k_2-2k_1+1}|q^2) \\ =\ & p_p(\cos\psi; qe^{i\te}/s, qe^{-i\te}/s,
q^{2k_2-2k_1+1}, q^{2k_1+2k_2-1} |q^2)\\ & \times
\ph{2}{2}{qe^{i(\te-\psi)}/s, qe^{i(\te+\psi)}/s } {q^2/s^2, q^{2k_2-2k_1-2m+2}e^{i\te}/s} {q^2,
q^{2k_2-2k_1-2m+2}e^{-i\te}/s },
\end{split}
\end{equation}
where
\[
\begin{split}
c_n&= \frac{(-1)^{n+m} q^{n(2k_2+1)} q^{n(n-1)-m(m-1)+2p(p-1)}} {
(q^2, q^{4k_2},q^{4k_1+2m}, q^{2+2m} ;q^2)_n}
(sq^{2k_1-2k_2-2})^{-m} (s^2q^{2k_1-2k_2-3})^{-p}  \\
&\times  \frac{ (q^{2+2m},q^{4k_1+2m};q^2)_\infty \,
(q^{4k_2};q^2)_p } { (q^2/s^2, q^{2k_2-2k_1-2m+2}
e^{i\te}/s;q^2)_\infty} \left| \frac{ (qe^{i(\psi+\te)}/s,
qe^{i(\psi-\te)}/s;q^2)_\infty } { (q^{2k_1+2k_2+2p}e^{i\te}/s,
q^{2k_1-2k_2+2m+1}e^{i\psi};q^2)_\infty } \right|^2.
\end{split}
\]
\end{thm}
Here we use the convention $s_{-n}=0$ for $n \geq 1$. So for $m <0$ the summation starts at $n=-m$. In this case the continuous dual $q$-Hahn polynomial $P_n$ is in fact a multiple of the polynomial $P_{n+m}(\cos\psi;q^{2k_2-2k_1-2m+1},q^{2k_1+2k_2-1},
q^{2k_1-2k_2+1}|q^2)$.

\begin{proof}
First we assume $m \geq 0$. Note that the $_2\varphi_2$-series in \eqref{sum1} can be written as a little $q$-Jacobi function $\phi_{m-L}(\cos\te;q^2/s^2, qe^{i\psi}/s, q^{2k_2-2k_1-2L+1}e^{-i\psi}|q^2)$, see remark \ref{rem1}.  We expand the product of the Askey-Wilson polynomial $p_p(\cos\psi;\cos\te)$ and the $_2\varphi_2$-series on the right hand side of \eqref{sum1} in terms of continuous dual $q$-Hahn polynomials $P_n(\cos\psi;q^{2k_1-2k_2+2m+1},q^{2k_1+2k_2-1},
q^{2k_2-2k_1+1}|q^2)$;
\[
p_p(\cos\psi;\cos\te)\phi_{m-L}(\cos\te;\cos\psi) = \sum_{n=0}^\infty C_n P_n(\cos\psi).
\]
Now we write the polynomials $p_p$ and $P_n$ in orthonormal form. Multiplying for fixed $k \in \Z_{\geq 0}$ both sides with $\Pt_k(\cos\psi;q^{2k_1-2k_2+2m+1},q^{2k_1+2k_2-1},
q^{2k_2-2k_1+1}|q^2)$ and integrating against the orthogonality measure $dm(\cos\psi;q^{2k_1-2k_2+2m+1},q^{2k_1+2k_2-1},
q^{2k_2-2k_1+1}|q^2)$, we find from the orthogonality of the continuous dual $q$-Hahn polynomials and proposition \ref{prop5}
\[
C_n = C F_{n+m,n}(\cos\te,p),
\]
where $C=C(\cos\psi, \cos\te, p, s, k_1, k_2)$ comes from the orthogonality measures. So now we have \eqref{sum1} with all the polynomials in orthonormal form. Writing all the polynomials in the usual normalization, we find after careful bookkeeping the desired summation formula for $m\geq 0$.

Next we assume $m \leq 0$. Applying the same method as above with the continuous dual $q$-Hahn polynomial $P_n(\cos\psi; q^{2k_2-2k_1-2m+1}, q^{2k_1+2k_2-1},  q^{2k_1-2k_2+1}|q^2)$, we find
\[
\begin{split}
\sum_{n=0}^\infty d_n s_{n}&(\cos\te;sq^{2k_1-2k_2-2p},
q^{2k_1+2k+2p}/s|q^2) s_{n+m}(\mu(q^{2k_2+2p}/s);sq^{-2k_2},q^{-2k_2}/s
|q^{-2}) \\ &  \times P_n(\cos\psi; q^{2k_2-2k_1-2m+1},q^{2k_1+2k_2-1},
q^{2k_1-2k_2+1}|q^2)  \\=\ & p_p(\cos\psi; qe^{i\te}/s, qe^{-i\te}/s,
q^{2k_2-2k_1+1}, q^{2k_1+2k_2-1} |q^2) \\ & \times
\ph{2}{2}{qe^{i(\te-\psi)}/s, qe^{i(\te+\psi)}/s } {q^2/s^2, q^{2k_2-2k_1-2m+2}e^{i\te}/s} {q^2,
q^{2k_2-2k_1-2m+2}e^{-i\te}/s },
\end{split}
\]
where
\[
\begin{split}
d_n&= \frac{(-1)^{n+m} q^{(n-m)(2k_2+1)}
q^{(n-m)(n-m-1)+2p(p-1)}} { (q^2, q^{4k_2-2m},q^{4k_1}, q^{2-2m}
;q^2)_n}
(q/s)^{m} (s^2q^{2k_1-2k_2-3})^{-p}  \\
&\times  \frac{ (q^{2-2m},q^{4k_1}, q^{4k_2-2m};q^2)_\infty \,
(q^{4k_2};q^2)_p } { (q^{4k_2}, q^2/s^2, q^{2k_2-2k_1-2m+2}
e^{i\te}/s;q^2)_\infty} \left| \frac{ (qe^{i(\psi+\te)}/s,
qe^{i(\psi-\te)}/s;q^2)_\infty } { (q^{2k_1+2k_2+2p}e^{i\te}/s,
q^{2k_1-2k_2+1}e^{i\psi};q^2)_\infty } \right|^2.
\end{split}
\]
Now we use \eqref{CDH} and the symmetry in the parameters $a$ and
$c$ to rewrite the continuous dual $q$-Hahn polynomial, cf. \eqref{eq:shift};
\[
\begin{split}
&\frac{ (q^{2-2m}, q^{4k_1};q^2)_\infty } { (q^{2-2m}, q^{4k_1};q^2)_n
} P_n(\cos\psi; q^{2k_2-2k_1-2m+1}, q^{2k_1+2k_2-1}, q^{2k_1-2k_2+1}
|q^2)= \\ &  (-1)^m q^{m(m-1)+2m}
(q^{2k_1-2k_2+2m+1})^{-m} \frac{ (q^2; q^2)_n (q^{4k_1+2m}, q^{2+2m}
;q^2)_\infty } { (q^2, q^{4k_1+2m}, q^{2+2m};q^2)_{n-m} } \\&
\times|(q^{2k_1-2k_2+1}e^{i\psi};q^2)_m |^2 P_{n-m}(\cos \psi;
q^{2k_1-2k_2+2m+1}, q^{2k_1+2k_2-1}, q^{2k_2-2k_1+1}|q^2).
\end{split}
\]
We now shift the summation index $n \mapsto n+m$ to see that
the expression for $m \leq 0$ is the same as for $m \geq 0$.
\end{proof}

\begin{thm} \label{thm2}
Let $m \in \Z$, $n \in \Z_{\geq 0}$ and $s \in \R$, with $|s| \geq q^{-1}$. For the Askey-Wilson polynomials $p_p$, the continuous dual $q$-Hahn polynomials $P_n$ and the Al-Salam and Chihara polynomials $s_n$, the following summation formula holds
\[
\begin{split}
\sum_{p=0}^\infty c_p s_{n+m}&(\cos\te;sq^{2k_1-2k_2-2p},
q^{2k_1+2k_2+2p}/s|q^2) s_{n}(\mu(q^{2k_2+2p}/s);sq^{-2k_2},q^{-2k_2}/s
|q^{-2}) \\  & \times  p_p(\cos\psi; qe^{i\te}/s, qe^{-i\te}/s,
q^{2k_2-2k_1+1}, q^{2k_1+2k_2-1} |q^2) \\=\ &P_n(\cos\psi;
q^{2k_1-2k_2+2m+1},q^{2k_1+2k_2-1}, q^{2k_2-2k_1+1}|q^2) \\ & \times
\ph{2}{2}{qe^{i(\te-\psi)}/s, qe^{i(\te+\psi)}/s } {q^2/s^2, q^{2k_2-2k_1-2m+2}e^{i\te}/s} {q^2,
q^{2k_2-2k_1-2m+2}e^{-i\te}/s },
\end{split}
\]
where
\[
\begin{split}
c_p =& \frac{ (q^{4k_2+2}/s^2;q^2)_{2p}
|(q^{2k_2-2k_1+2p+2}e^{i\te}/s)_\infty |^2 } {(q^2, q^2/s^2,
q^{4k_2+2p}/s^2 ;q^2)_p } ( q^{2k_2-2k_1-3}/s^2 )^p q^{2p(p-1)} \\
&\times \frac{ (-1)^{n+m} q^{n(1+2k_2)} (sq^{2k_1-2k_2})^{-m}
q^{n(n-1)-m(m-1)} } { ( q^{4k_2+2}/s^2, q^{2k_1-2k_1-2m+2}e^{i\te}/s
;q^2)_\infty }.
\end{split}
\]
\end{thm}
Again the convention $s_{-n}=0$, $n \geq 1$, is used.
\begin{proof}
The proof runs along the same lines as the proof of theorem \ref{thm1}.
\end{proof}
\begin{rem}
(i) In view of remark \ref{Jac pol} we can replace the Al-Salam polynomial in base $q^{-2}$ in theorems \ref{thm1} and \ref{thm2} by a little $q$-Jacobi polynomial in base $q^2$;
\[
\begin{split}
s_{n}(&\mu(q^{2k_2+2p}/s);sq^{-2k_2},q^{-2k_2}/s
|q^{-2}) = \\ & (-1)^{n+p} q^{-n(n-1)-p(p+1)} s^{2p} q^{-2nk_2} (s^2 q^{4k_2};q^2)_{n-p}\, (q^2/s^2;q^2)_p \,p_p(q^{2n}; s^{-2}, q^{4k_2-2};q^2).
\end{split}
\]

(ii) An expression also involving Al-Salam and Chihara polynomials, with a structure that is similar to, but simpler than, the expressions in theorems \ref{thm1} and \ref{thm2}, can be found in \cite[Thm.4.3]{IS}. The expression in \cite{IS} does not seem to be related to the expressions in this paper, since it involves two Al-Salam and Chihara polynomials in base $q$.

(iii) Theorems \ref{thm1} and \ref{thm2} may both be considered as
$q$-analogues of \cite[Theorem 3.6]{GK}.

(iv) Theorem \ref{thm1} gives the inverse connection formula for the two orthogonal bases mentioned in remark \ref{rem2}(iii) on $\mathcal D$, see also figure \ref{fig1}.

(v) If we take $n=m=0$ in theorem \ref{thm2}, we find the
summation formula \eqref{expan}.
\end{rem}

\section{Holomorphic and anti-holomorphic realizations} \label{sec4}
In this section we realize the basisvectors for the representation spaces of the discrete series representations as holomorphic and anti-holomorphic functions. In these realizations the eigenvectors $v^+$ and $v^-$, \eqref{v+} and \eqref{v-}, become known generating functions. Using the realizations of the standard basisvectors and of the eigenvectors, we find a bilinear generating function for a special type of $_2\varphi_1$-series. See also \cite{JJ} and \cite{Ros} where similar realizations are being used to find generating functions related to positive discrete series representations.

To simplify notations we assume $k_1-k_2 \geq -\frac{1}{2}$ and $k_1+k_2 \geq
\frac{1}{2}$, then the measure $dm^\hf$ in the Clebsch-Gordan decomposition \eqref{CG decomp} has only an absolutely continuous part. Further we assume without loss of generality $x=\mu(u)$ with $|u|\leq 1$, where $(2x-2\mu(s))/(q^{-1}-q)$ is a spectral value of $\pi^+_{k_1}\tensor \pi^-_{k_2}\big(\De(Y_sA)\big)$, cf. proposition \ref{prop4}. In case discrete terms occur in the decomposition \eqref{CG decomp}, the calculations proceed completely analogous.

\subsection{Holomorphic and anti-holomorphic realizations}
Consider the Hilbert space $\mathcal{H}_q^k$ of holomorphic functions on the unit disk $\{z \in \C ; \,|z|<1\}$ with finite norm with respect to the inner product
\[
\inprod{f}{g}= \sum_{n=0}^\infty \frac{ (q^2;q^2)_n }{(q^{4k};q^2)_n}
f_n \overline{g}_n, \qquad f(z)= \sum_{n=0}^\infty f_n z^n.
\]
Standard orthonormal basisvectors in this space are
\begin{equation} \label{rea1}
e_n =  \sqrt{ \frac{(q^{4k};q^2)_n}{ (q^2;q^2)_n } } z^n.
\end{equation}
The realization of the positive discrete series representation on
the space $\mathcal{H}_q^k$ can be given in terms of the
dilatation operator $T_q$ and $q$-derivative operator $D_q$
given by
\[
T_q f(z) = f(qz), \qquad D_q = \frac{1-T_q}{(1-q)z}.
\]
The realization is
\[
\begin{split}
\pi^+_k (A) &= q^k T_q, \quad \pi^+_k(D) = q^{-k} T_{q^{-1}},\\
\pi^+_k(B) & = q^{\hf(1-2k)}\left(z^2D_{q^2}T_{q^{-1}} + \frac{
1-q^{4k}}{1-q^2} z T_q\right), \\ \pi^+_k(C) &=-q^{\hf(3-2k)}
D_{q^2} T_{q^{-1}}.
\end{split}
\]

Consider the space $\overline{\mathcal{H}}_q^k$ of anti-holomorphic functions on the unit disk with finite norm with respect to the inner product
\[
\inprod{f}{g}= \sum_{n=0}^\infty \frac{ (q^{-2};q^{-2})_n
}{(q^{-4k};q^{-2})_n} f_n \overline{g}_n, \qquad f(z)=
\sum_{n=0}^\infty f_n \overline{w}^n.
\]
Standard orthonormal basisvectors in this space are
\begin{equation} \label{rea2}
e_n = \sqrt{ \frac{(q^{-4k};q^{-2})_n }{ (q^{-2};q^{-2})_n} }
\overline{w}^n.
\end{equation}
The realization of the negative discrete series on the space
$\overline{\mathcal{H}}_q^k$ is
\[
\begin{split}
\pi^-_k(A) & = q^{-k}T_{q^{-1}},\qquad \pi^-_k(D)=q^kT_q ,\\
\pi^-_k(B) &= q^{\hf(1+2k)}D_{q^2}T_{q^{-1}},\\ \pi^-_k(C)&=
-q^{\hf(3+2k)}\left(\w^2D_{q^2}T_{q^{-1}} + \frac{
1-q^{-4k}}{1-q^2} \w T_{q^{-1}}\right).
\end{split}
\]

Next we use the realization of the standard orthonormal
basisvectors \eqref{rea1} and \eqref{rea2} to find an expression
for $f \tensor e_n \in L^2(-1,1)\tensor \ell^2(\Z)$ as a function of $z$ and $\w$.
\begin{prop} \label{propo1}
For $n \in \Z$ and $|z\w|<q^{2k_2-1}$
\[
\begin{split}
(f\ \tensor&\ e_{n-L})(z,\w) = (-1)^n q^{n(n-1)} (\w q^{-2k_1})^{-n}
(z\w q^{2k_1}, q^{2-2n};q^2)_\infty
\sqrt{\frac{(q^{4k_1},q^{4k_2};q^2)_\infty}{2\pi}} \\ \times&
\int_0^\pi \frac{f(\cos\psi)}{(z\w q^{1-2k_2}e^{-i\psi};q^2)_\infty}
\ph{2}{1}{q^{2k_2-2k_1-2n+1}e^{-i\psi}, q^{2k_1-2k_2+1}e^{-i\psi}}
{q^{2-2n}}{q^2,z\w q^{1-2k_2}e^{i\psi}} \\& \times    \left| \frac{
(e^{2i\psi};q^2)_\infty } {(q^{2k_1-2k_2+1}e^{i\psi};q^2)_n
(q^{2k_2-2k_1+1}e^{i\psi}, q^{2k_1+2k_2-1}e^{i\psi},
q^{2k_1-2k_2+1}e^{i\psi};q^2)_\infty} \right| d\psi.
\end{split}
\]
\end{prop}
\begin{proof}
We start with the case $n \leq 0$. We substitute \eqref{rea1} and \eqref{rea2} in the decomposition
\eqref{inv CG} and interchange summation and integration. This is
allowed for a sufficiently smooth function $f$. Next we write out the
measure $dm^\hf$ explicitly and we write the continuous dual $q$-Hahn
polynomial in the normalization given by \eqref{CDH}. Then we use the
generating function for continuous dual $q$-Hahn polynomials, see \cite{KS},
\[
\sum_{m=0}^\infty \frac{ P_m(\cos\psi;a,b,c|q) }{ (q,ac;q)_m} t^m = \frac{ (bt;q)_\infty }{ (te^{-i\psi};q)_\infty } \ph{2}{1} {ae^{-i\psi}, ce^{-i\psi} }{ac}{q,te^{i\psi}}, \qquad |t|<1,
\]
to evaluate the sum. Note that this generating function can be derived directly from theorem \ref{thm2.1} by first putting $b=0$ in theorem \ref{thm2.1} and then $r=0$. Finally we rewrite the one term inside the modulus
signs that depends on $n$, using
\begin{equation} \label{sh}
(aq^{-n};q)_n = (q/a;q)_n \left(-\frac{a}{q}\right)^n q^{-\hf n(n-1)}.
\end{equation}
This gives the desired expression for $n \leq 0$.

For $n\geq 0$ we apply the same method, apart from the last step, to
find
\[
\begin{split}
(f\ \tensor&\ e_{n-L})(z,\w) =  z^n (z\w q^{2k_1}, q^{2+2n};q^2)_\infty
\sqrt{\frac{(q^{4k_1},q^{4k_2};q^2)_\infty}{2\pi}}
\\ \times& \int_0^\pi \frac{f(\cos\psi)}{(z\w
q^{1-2k_2}e^{i\psi};q^2)_\infty} \ph{2}{1}{q^{2k_1-2k_2+2n+1}e^{-i\psi},
q^{2k_2-2k_1+1}e^{-i\psi}} {q^{2+2n}}{q^2,z\w q^{1-2k_2}e^{i\psi}} \\&
\times    \left| \frac{(q^{2k_1-2k_2+1}e^{i\psi};q^2)_n
(e^{2i\psi};q^2)_\infty } {(q^{2k_2-2k_1+1}e^{i\psi},
q^{2k_1+2k_2-1}e^{i\psi}, q^{2k_1-2k_2+1}e^{i\psi};q^2)_\infty} \right|
d\psi.
\end{split}
\]
We use the following transformation to rewrite the
$_2\varphi_1$-series; for $n \in \Z$
\begin{equation} \label{eqB}
(q^{1+n};q)_\infty \ph{2}{1}{aq^n,b}{q^{1+n}}{q,t} = (-1)^n q^{\hf
n(n-1)}\frac{(q^{1-n};q)_\infty}{(a,q/b;q)_n}
\left(\frac{q}{bt}\right)^n \ph{2}{1}{a, bq^{-n}}{q^{1-n}}{q,t}.
\end{equation}
This behaviour of a $_2\varphi_1$-series is usually called Bessel
coefficient behaviour. Finally using \eqref{sh} on
$(q^{2k_2-2k_1-2n+1};q^2)_n$ and cancelling common factors gives
the same expression as for $n\leq 0$.
\end{proof}

\subsection{Realizations of the eigenvectors of $\boldsymbol{Y_sA}$}
Using the orthonormal basisvectors \eqref{rea1}, the generalized
eigenvector $v^+$, see \eqref{v+}, becomes a generating function for the
orthonormal Al-Salam and Chihara polynomials (see
\cite[eq.(3.10)]{AI}, \cite{KS})
\begin{equation} \label{rea v+}
v^+(\cos\te;z)=
\sum_{n=0}^\infty \sqrt{ \frac{(q^{4k};q^2)_n}{ (q^2;q^2)_n } }
\s_n(\cos\te;sq^{2k}, q^{2k_1}/s|q^2) z^n = \frac{(szq^{2k},
zq^{2k}/s;q^2)_\infty} { (ze^{i\te}, ze^{-i\te};q^2)_\infty}.
\end{equation}
Using \eqref{rea2} the eigenvector $v^-$, see \eqref{v-}, also becomes a generating function (see \cite[eq.(3.70)]{AI})
\begin{equation} \label{rea v-}
\begin{split}
v^-(p;\w) &= \sum_{n=0}^\infty (-1)^n \sqrt{ \frac{(q^{-4k};q^{-2})_n }{
(q^{-2};q^{-2})_n} }
\s_n\big(\mu(q^{2k+2p}/s);sq^{-2k},q^{-2k}/s|q^{-2}\big) \w^n
\\&= \frac{ (s\w q^{2-2k-2p}, \w q^{2k+2p+2}/s
;q^2)_\infty} { (s\w q^{2-2k}, \w q^{2-2k}/s
;q^2)_\infty}.
\end{split}
\end{equation}
These generating functions enable us to give an explicit expression for the uncoupled eigenvector $v(x,p)$. From \eqref{v(xp)} and the explicit expressions for $v^+$ and $v^-$, \eqref{rea v+} and \eqref{rea v-}, we find for $x=\mu(u)$
\begin{equation} \label{expl v}
\begin{split}
v&(x,p;z,\w) = \frac{
(szq^{2k_1-2k_2-2p}, zq^{2k_1+2k_2+2p}/s, s\w q^{2-2k_2-2p}, \w
q^{2k_2+2p+2}/s; q^2)_\infty } { (zu, z/u, s\w q^{2-2k_2} , \w
q^{2-2k_2}/s ;q^2)_\infty } \\
& =\frac{ (szq^{2k_1-2k_2}, zq^{2k_1+2k_2}/s, \w q^{2k_2+2}/s ;q^2)_\infty \, (q^{2k_2-2k_1+2}/sz, q^{2k_2}/s\w ;q^2)_p } { (zu, z/u, \w q^{2-2k_2}/s ;q^2)_\infty \, (zq^{2k_1+2k_2}/s, \w q^{2k_2+2}/s ;q^2)_p }\\ & \qquad\times (s^2z\w q^{2k_1-4k_2-2})^p q^{-2p(p-1)}.
\end{split}
\end{equation}
Here we use
\[
(aq^p;q)_\infty = \frac{ (a;q)_\infty }{ (a;q)_p }, \qquad (aq^{-p};q)_\infty = (-a/q)^p q^{-\hf p(p-1)} (q/a;q)_p \, (a;q)_\infty, \qquad p \in \Z_{\geq 0}.
\]

Next we use proposition \ref{propo1} to give an explicit expression for the coupled eigenvectors
\[
f \tensor v^P(x) = \int_0^\pi f(\cos\psi) v^P(x) d\psi
\]
in this realization.
\begin{prop} \label{propo1a}
For $|\w q^{2-2k_2}/s|<1$ and $x = \mu(u)$
\begin{equation} \label{fv}
\begin{split}
\big(f \tensor& v^P(x) \big) (z,\w)=  (-q/s)^L (z\w q^{2k_1}, z\w q^{2k_1-4k_2+2};q^2)_\infty   \sqrt{ \frac{ (q^{4k_1}, q^{4k_2};q^2)_\infty }{ 2\pi }} \\
&\times \int_0^\pi \frac{f(\cos \psi)}{(z\w
q^{1-2k_2}e^{i\psi},z\w q^{1-2k_2}e^{-i\psi};q^2)_\infty } \left| \frac{ (e^{2i\psi};q^2) } { (q^{2k_1+2k_2-1}e^{i\psi}, q^{2k_1-2k_2+1}e^{i\psi} ;q^2)_\infty } \right| \\
&\times  \sum_{n=-\infty}^\infty  \ph{2}{1}{z\w q^{1-2k_2}e^{-i\psi}, q^{2k_1-2k_2+1}e^{-i\psi}} {z\w q^{2k_1-4k_2+2}}{q^2,q^{2k_2-2k_1-2n+1} e^{i\psi}}\\
&\qquad  \times  \ph{2}{1}{que^{i\psi}/s, qe^{i\psi}/us} {q^2/s^2}{q^2, q^{2k_2-2k_1-2n+1}e^{-i\psi}} (sq^{2k_2-2}/\w)^n\,d\psi.
\end{split}
\end{equation}
\end{prop}
\begin{proof}
From \eqref{vP} we find for $x=\mu(u)$
\[
\big(f \tensor v^P(x) \big) (z,\w)= \int_0^\pi f(\cos \psi)
\left(\sum_{n=-\infty}^\infty \phi_{n-L}(x;q^2/s^2, qe^{i\psi}/s,
q^{2k_2-2k_1+1}e^{-i\psi}|q^2) e_{n-L}\right) (z,\w)d\psi.
\]
Then \eqref{def Jac} and proposition \ref{propo1} give the desired expression for $(f \tensor v^P)(z,\w)$ after using Heine's transformation \cite[eq.(III.2)]{GR} and \eqref{sh}.

To determine conditions for absolute convergence of the inner sum in
\eqref{fv}, we investigate the asymptotic behaviour of the summand. For $n\rightarrow -\infty$ we see that the summand is $\mathcal
O\big((sq^{2k_2-2}/\w)^n\big)$. So we need $|\w q^{2-2k_2}/s|<1$ for absolute convergence. For $n \rightarrow \infty$ we use the connection formula \cite[eq.(A.10)]{KSt2} to find that the first $_2\varphi_1$-series behaves like
$A(que^{i\psi}/s)^n + B(qe^{i\psi}/us)^n$, where $A$ and $B$ are nonzero constants for $x=\mu(u) \in [-1,1]$ and $A\neq 0$, $B=0$ for $u<1$, assuming that $x \in \mathrm{supp} \big( d\nu(\cdot;q^2/s^2, qe^{i\psi}/s, q^{2k_2-2k_1-2L+1}e^{-i\psi}|q^2)\big)$. We transform the second $_2\varphi_1$-series back to the form it has in proposition \ref{propo1} using \cite[eq.(III.2)]{GR}. Then we use the Bessel coefficient behaviour \eqref{eqB} and transformation \cite[eq.(III.2)]{GR} again, to find
\[
\begin{split}
( q^{2k_2-2k_1-2n+1}&e^{i\psi}, z\w q^{2k_1-4k_2+2} ;q^2)_\infty  \ph{2}{1}{z\w q^{1-2k_2}e^{-i\psi}, q^{2k_1-2k_2+1}e^{-i\psi}} {z\w q^{2k_1-4k_2+2}}{q^2,q^{2k_2-2k_1-2n+1} e^{i\psi}} = \\&
(-1)^n q^{-n(n-1)} (z\w q^{-2k_1})^n |  (q^{2k_1-2k_2+1}e^{i\psi};q^2)_n |^2  (q^{2k_1-2k_2+2n+1}e^{i\psi}, z\w q^{2-2k_1};q^2)_\infty \\
&\times \ph{2}{1} {z\w q^{1-2k_2}e^{-i\psi}, q^{2k_2-2k_1+1}e^{-i\psi}} {z\w q^{2-2k_1}} {q^2, q^{2k_1-2k_2+2n+1}e^{i\psi} }.
\end{split}
\]
We use
\[
(a;q)_n = (q^{1-n}/a;q)_n (-a)^n q^{\hf n(n-1)}
\]
to see that summand is $\mathcal O((uz)^n)$, $n \rightarrow \infty$. Recall that we assume $|u|\leq 1$ and $|z|<1$, so for $n \rightarrow \infty$ the conditions for absolute convergence are satisfied.
\end{proof}

Next we determine another expression for $(f \tensor v^P)(z,\w)$, using the Clebsch-Gordan decomposition for the eigenvectors \eqref{CG ev}. The method we use is in principle the same as for the standard orthonormal basisvectors: we determine $f \tensor v^P(x)$ from the explicit realization for $v(x,p)$, see \eqref{v(xp)}, and the inverse formula to \eqref{CG ev}.
Define a new measure $d\M(y;x,p)$ by
\[
d\M(y;x,p) d\nu(x; q^2/s^2, qe^{i\psi}/s, q^{2k_2-2k_1-2L+1}e^{-i\psi}|q^2)= dm(y;a,b,c,d) dm(x,p), \quad y=\cos\psi,
\]
where $a$, $b$, $c$, $d$ are given by \eqref{abcd1}.
Inverting \eqref{CG ev} we find for sufficiently smooth $f$
\begin{equation} \label{fv inv}
f\, \tensor \,v^P(x)  =(-\sgn s)^L \sum_{p=0}^\infty \left[ \int_\R \p_p(y;x) f(y) d\M^\hf(y;x,p)\right]  v(x,p) .
\end{equation}
This expression enables us to evaluate $(f \tensor v^P)(z,\w)$ in a different way.
\begin{prop} \label{prop4.2}
For $x=\mu(u)$ and $|z\w | < q^{2k_2-1}$
\[
\begin{split}
\big(f \tensor& v^P(x) \big) (z,\w) =
\frac{ (q^2, z\w q^{2k_1}, szq^{2k_1-2k_2}, uq^{2k_2-2k_1+2}/s, q^{2k_2-2k_1+2}/su ;q^2)_\infty } { (\w q^{2-2k_2}/s, zu,z/u;q^2)_\infty } \\
\times & (-q/s)^L \sqrt{\frac{(q^{4k_1}, q^{4k_2};q^2)_\infty } {2\pi}} \int_0^\pi  f(\cos\psi) \Phi(\cos\psi) \left| \frac{ (e^{2i\psi};q^2)_\infty } { (q^{2k_1+2k_2-1}e^{i\psi}, q^{2k_1-2k_2+1}e^{i\psi} ;q^2)_\infty } \right|  d\psi,
\end{split}
\]
where
\[
\begin{split}
\Phi&(\cos\psi) = \frac{(\w q^{3-2k_1}e^{i\psi}/s, zqe^{i\psi}/s ;q^2)_\infty} { (q^{2k_2-2k_1+1}e^{i\psi}, q^{2k_2-2k_1+1}e^{-i\psi},  q^{2k_2-2k_1+3}e^{i\psi}/s^2, z\w q^{1-2k_2}e^{i\psi} ;q^2)_\infty }  \\\times&_8W_7(q^{2k_2-2k_1+1}e^{i\psi}/s^2; que^{i\psi}/s, qe^{i\psi}/su, q^{2k_2-2k_1+1}e^{i\psi}, q^{2k_2}/s\w, q^{2k_2-2k_1+2}/sz;q^2, z\w q^{1-2k_2}e^{-i\psi}).
\end{split}
\]
\end{prop}
\begin{proof}
We insert the explicit expression for the uncoupled eigenvector $v(x,p;z,\w)$
\eqref{expl v} in \eqref{fv inv} and we interchange summation and
integration. This is allowed if we choose $f$ to be sufficiently
smooth. Next we write the Askey-Wilson polynomial in the usual
normalization and we bring all terms not depending on $p$ outside the
sum. After some calculations we find
\[
\begin{split}
\big(f \tensor &v^P(x) \big) (z,\w) = (-q/s)^L \sqrt{\frac{(q^{4k_1}, q^{4k_2};q^2)_\infty } {2\pi}}\\
\times & \frac{ (q^2, szq^{2k_1-2k_2}, zq^{2k_1+2k_2}/s, \w q^{2k_2+2}/s;q^2)_\infty } { (\w q^{2-2k_2}/s, q^{4k_2+2}/s^2;q^2)_\infty } \frac{ | (q^{2k_2-2k_1+2}e^{i\te}/s;q^2)_\infty |^2} { (ze^{i\te},ze^{-i\te};q^2)_\infty }  \\
\times & \int_0^\pi \frac{ f(\cos\psi) } { |(q^{2k_2-2k_1+1}e^{i\psi};q^2)_\infty |^2 }\left| \frac{ (e^{2i\psi};q^2)_\infty } { (q^{2k_1+2k_2-1}e^{i\psi}, q^{2k_1-2k_2+1}e^{i\psi} ;q^2)_\infty } \right| S(\cos\psi) d\psi,
\end{split}
\]
where $S(\cos\psi)= S(\cos\psi, u, s, k_1, k_2, z, \w)$ is the sum over $p$;
\[
\begin{split}
S(\cos\psi) = \sum_{p=0}^\infty &\frac{ (q^{4k_2+2}/s^2;q^2)_{2p}\, \, p_p(\cos\psi;qu/s, q/su, q^{2k_2-2k_1-1}, q^{2k_1+2k_2-1}|q^2) } { (q^2, q^2/s^2, uq^{2k_2-2k_1+2}/s, q^{2k_2-2k_1+2}/su, q^{4k_2+2p}/s^2 ;q^2)_p } \\
& \times \frac{ (q^{2k_2-2k_1+2}/sz, q^{2k_2}/s\w ;q^2)_p }{ (zq^{2k_1+2k_2}/s, \w q^{2k_2+2}/s;q^2)_p } (z\w q^{1-2k_2})^p.
\end{split}
\]
Here we recognize the generating function from theorem \ref{thm2.1} with $a$, $b$, $c$, $d$ as in \eqref{abcd1} and with $r = zq/s$, $t= z\w q^{1-2k_2}$.
So this gives
\[
\begin{split}
S&(\cos\psi)= \frac{ (q^{4k_2+2}/s^2, z\w q^{2k_1}, \w q^{3-3k_1}e^{i\psi}/s, zqe^{i\psi}/s ;q^2)_\infty } { (\w q^{2k_2+2}/s, zq^{2k_1+2k_2}/s, q^{2k_2-2k_1+3}e^{i\psi}/s^2, z\w q^{1-2k_2}e^{i\psi} ;q^2)_\infty }\\  \times
{}_8W_7&(q^{2k_2-2k_1+1}e^{i\psi}/s^2; que^{i\psi}/s, qe^{i\psi}/su, q^{2k_2-2k_1+1}e^{i\psi}, q^{2k_2}/s\w, q^{2k_2-2k_1+2}/sz;q^2, z\w q^{1-2k_2}e^{-i\psi}).
\end{split}
\]
This proves the proposition.
\end{proof}

\subsection{A bilinear generating function}
Propositions \ref{propo1a}  and \ref{prop4.2} both give an expression for $(f \tensor v^P)(z,\w)$. Since $f$ is arbitrary, these two expression must be equal. This gives a bilinear summation formula.
\begin{thm} \label{thm4.4}
For $|a|\leq 1$, $b<1$, $|u|\leq 1$ and $1<|t|<|abu|^{-1}$
\[
\begin{split}
&\sum_{n=-\infty}^\infty  \ph{2}{1}{be^{-i\psi},  qe^{-i\psi}/c}{bq/c}{q,cq^{-n}e^{i\psi}}\ph{2}{1}{ae^{i\psi}/u, aue^{i\psi}}{a^2}{q,cq^{-n}e^{-i\psi}}  t^n=\\&
\frac{ (q, be^{-i\psi}, ce^{i\psi}/t, a^2bte^{i\psi}, bct, acu, ac/u ;q)_\infty }{ (ce^{i\psi}, ce^{-i\psi}, a^2ce^{i\psi}, bq/c, 1/t, abtu, abt/u ;q)_\infty}\,
{}_8W_7(a^2ce^{i\psi}/q; aue^{i\psi}, ae^{i\psi}/u, ce^{i\psi}, a^2t, c/bt;q, be^{-i\psi})
\end{split}
\]
\end{thm}
\begin{proof}
This follows from propositions \ref{propo1a} and \ref{prop4.2}, relabelling
\[
a \mapsto q/s, \qquad
b \mapsto z\w q^{1-2k_2},\qquad c \mapsto q^{2k_2-2k_1+1}, \qquad t \mapsto sq^{2k_2-2}/\w
\]
and replacing $q^2$ by $q$.
\end{proof}
\begin{rem}
The expression in theorem \ref{thm4.4} is first proved by Rahman in \cite[App.~B.3]{KSt2} by analytic methods. We can consider theorem \ref{thm4.4} as a special case of the dual transmutation kernel for the little $q$-Jacobi functions, see \cite[Thm.~2.1]{KR}. The general dual transmutation kernel is expressed as the sum of two very-well-poised $_8\varphi_7$-series. In the case of theorem \ref{thm4.4} one of the $_8\varphi_7$-series vanishes. This can be easily seen using \cite[eq.(3.12)]{KR}.
\end{rem}

\end{document}